\documentclass[9pt]{article}
\usepackage{mathrsfs}
\usepackage{amsthm}
\usepackage{amssymb}
\usepackage{amsmath}
\usepackage{graphicx}
\usepackage{color}
\usepackage{amsfonts}
\usepackage{float}
\usepackage{cite}
\usepackage [latin1]{inputenc}
\usepackage[text={140mm,210mm},left=35mm,vmarginratio=1:1]{geometry}
\newtheorem{theorem}{Theorem}[section]
\newtheorem{remark}{Remark}[section]
\newtheorem{lemma}[theorem]{Lemma}

\numberwithin{equation}{section}
\normalsize

\begin{document}
\title{\textbf{Central limit theorems and moderate deviations for additive functionals of SSEP on regular trees}}

\author{Xiaofeng Xue \thanks{\textbf{E-mail}: xfxue@bjtu.edu.cn \textbf{Address}: School of Mathematics and Statistics, Beijing Jiaotong University, Beijing 100044, China.}\\ Beijing Jiaotong University}

\date{}
\maketitle

\noindent {\bf Abstract:} In this paper, we are concerned with the symmetric simple exclusion process (SSEP) on the regular tree $\mathcal{T}_d$. A central limit theorem and a moderate deviation principle of the additive functional of the process are proved, which include the CLT and the MDP of the occupation time as special cases. A graphical representation of the SSEP plays the key role in proofs of the main results, by which we can extend the martingale decomposition formula introduced in Kipnis (1987) for the occupation time to the case of general additive functionals.

\quad

\noindent {\bf Keywords:} central limit theorem, moderate deviation, additive functional, exclusion process, graphical representation.

\section{Introduction}\label{section one}
This paper is a further investigation of the problem studied in \cite{Xue2025}. We extend the occupation time moderate deviation principle (MDP) of the symmetric simple exclusion process (SSEP) on the regular tree to the case of general additive functionals and a central limit theorem (CLT) is also proved. We first introduce some notations for later use. For each $d\geq 2$, we denote by $\mathcal{T}_d$ the regular tree where each vertex has $d+1$ neighbors. When $x, y\in \mathcal{T}_d$ are neighbors, we write $x\sim y$. According to the structure of $\mathcal{T}_d$, there exists $\alpha: \mathcal{T}_d\rightarrow \mathbb{Z}$ such that, for any $x\in \mathcal{T}_d$, there are $d$ neighbors $y$ of $x$ such that $\alpha(y)=\alpha(x)+1$ and there is one neighbor $z$ of $x$ such that $\alpha(z)=\alpha(x)-1$, i.e., each $x\in \mathcal{T}_d$ has $d$ children and one father. For any $x, y\in \mathcal{T}_d$, there is an unique integer $m\geq 0$ and an unique self-avoiding path $x=x_0\sim x_1\sim\ldots\sim x_m=y$ on $\mathcal{T}_d$. We denote $m$ by $D(x, y)$, i.e., $D(x, y)$ is the distance between $x$ and $y$.

Now we recall the definition of the SSEP. The SSEP $\{\eta_t\}_{t\geq 0}$ on $\mathcal{T}_d$ is a continuous-time Markov process with state space $\{0, 1\}^{\mathcal{T}_d}$. The generator $\mathcal{L}$ of $\{\eta_t\}_{t\geq 0}$ is given by
\begin{equation}\label{equ 1.1 generator}
\mathcal{L}f(\eta)=\frac{1}{2}\sum_{x\in \mathcal{T}_d}\sum_{y\sim x}\left(f(\eta^{x, y})-f(\eta)\right)
\end{equation}
for any local $f: \mathcal{T}_d\rightarrow \mathbb{R}$ and $\eta\in \{0, 1\}^{\mathcal{T}_d}$, where
\[
\eta^{x, y}(z)=
\begin{cases}
\eta(z) & \text{~if~}z\neq x, y,\\
\eta(x) & \text{~if~}z=y,\\
\eta(y) & \text{~if~}z=x.
\end{cases}
\]
According to the expression of $\mathcal{L}$, for a pair of neighbors $x$ and $y$, in the SSEP the value of $x$ and the value of $y$ are exchanged at rate $1$. If we consider $x$ in state $1$ as a vertex with a car on it and $y$ in state $0$ as a vacant  vertex, then all cars perform independent simple random walks on $\mathcal{T}_d$ except that each jump to an occupied vertex is suppressed. We can equivalently consider that all cars perform independent simple random walks except that any two cars on a pair of neighbors exchange their positions with each other at rate one. In section \ref{section three}, we recall a graphic representation of the SSEP, where the latter explanation is utilized. For a detailed survey of the basic properties of the SSEP, see Chapter 8 of \cite{Lig1985} and Part 3 of \cite{Lig1999}.

For any probability measure $\mu$ on $\{0, 1\}^{\mathcal{T}_d}$, we denote by $\mathbb{P}_\mu$ the probability measure of the SSEP $\{\eta_t\}_{t\geq 0}$ starting from $\mu$. We denote by $\mathbb{E}_\mu$ the expectation with respect to $\mathbb{P}_{\mu}$. For any $\eta\in \{0, 1\}^{\mathcal{T}_d}$, if $\mu$ is the Dirac measure concentrated on $\eta$, then we write $\mathbb{P}_\mu, \mathbb{E}_\mu$ as $\mathbb{P}_\eta, \mathbb{E}_\eta$ respectively. For any $p\in (0, 1)$, we denote by $\nu_p$ the product measure on $\mathcal{T}_d$ with density $p$, i.e., $\{\eta(x)\}_{x\in \mathcal{T}_d}$ are independent under $\nu_p$ and
\[
\nu_p\left(\eta(x)=1\right)=p=1-\nu_p(\eta(x)=0)
\]
for all $x\in \mathcal{T}_d$. According to the expression of the generator $\mathcal{L}$, it is easy to check that $\nu_p$ is a reversible distribution of $\{\eta_t\}_{t\geq 0}$.

In this paper, we are concerned with the additive functional of the SSEP. From now on, we let $p\in (0, 1)$ be a fixed constant and assume that $F: \{0, 1\}^{\mathcal{T}_d}\rightarrow \mathbb{R}$ be a fixed local function such that
\begin{equation}\label{equ basic mean zero assumption}
\int_{\{0, 1\}^{\mathcal{T}_d}} F(\eta) \nu_p(d\eta)=0.
\end{equation}
Let
\begin{equation}\label{equ 1.2 additive functional}
\xi_t^F=\int_0^tF(\eta_s)ds,
\end{equation}
which is called the `additive functional', then we have the following variance limit theorem.
\begin{lemma}\label{lemma 1.1}
Let $d\geq 2$. For any local function $F: \{0, 1\}^{\mathcal{T}_d}\rightarrow \mathbb{R}$ satisfying \eqref{equ basic mean zero assumption}, there exists $\sigma^2_F<+\infty$ such that
\[
\lim_{t\rightarrow+\infty}\frac{1}{t}{\rm Var}_{\nu_p}\left(\xi_t^F\right)=\sigma_F^2.
\]
\end{lemma}
The proof of Lemma \ref{lemma 1.1} is given in Section \ref{section three}, where the graphical representation of the SSEP introduced in \cite{Faggionato2024} is utilized. By Lemma \ref{lemma 1.1}, it is natural to ask whether a CLT and a MDP hold for $\xi_t^F$. In this paper, we give a positive answer to the above question. For a precise statement of our main results, see Section \ref{section two}.

A special and important case of the additive functional is the `occupation time', where $F(\eta)=\eta(x)-p$ for some $x\in \mathcal{T}_d$. The investigation of the occupation time limit theorems of the SSEP dates back to 1980s. Reference \cite{Kipnis1987}, by introducing a martingale decomposition strategy, gives the occupation time CLT of the SSEP on the lattice $\mathbb{Z}^d$. An interesting phenomenon in the $\mathbb{Z}^d$ case is that the scaling function of the CLT depends on the dimension $d$, which is $\sqrt{t}$ when $d\geq 3$, $\sqrt{t\log t}$ when $d=2$ and $t^{\frac{3}{4}}$ when $d=1$. This dimension-dependent phase transition relies on the fact that the simple random walk on $\mathbb{Z}^d$ is transient when $d\geq 3$ and recurrent when $d\leq 2$. By introducing a deviation inequality, Reference \cite{Gao2024} proves the occupation time MDP of the SSEP on $\mathbb{Z}^d$, where a dimension-dependent phase transition also occurs. Reference \cite{Xue2025} proves the occupation time MDP of the SSEP on $\mathcal{T}_d$, where the scaling function does not depend on $d$, which relies on the fact that the simple random walk on $\mathcal{T}_d$ is transient for all $d\geq 2$. Reference \cite{Landim1992} proves a process level large deviation principle (LDP) of the SSEP on $\mathbb{Z}^d$ and then, by utilizing the contraction principle, obtains the occupation time LDP in the cases of $d=1$ and $d\geq 3$. Reference \cite{Lee2004} further proves the occupation time LDP of the SSEP on $\mathbb{Z}^2$.

Limit theorems of general additive functionals of the SSEP are also discussed in previous literatures. By giving a sharp estimation on the `spectral gap' of the process, Reference \cite{Sethuraman1996} proves the CLT of the additive functionals of the SSEP on $\mathbb{Z}^d$. Reference \cite{Sethuraman2000} further extends the result in \cite{Sethuraman1996} to the asymmetric case. By utilizing a local average approach, Reference \cite{Gao2024} proves the MDP of the additive functionals of the SSEP on $\mathbb{Z}^1$.

\section{Main results}\label{section two}
In this section, we give our main results. Our first main result is about the CLT of the additive functional of $\{\eta_t\}_{t\geq 0}$. We let $T>0$ be a fixed moment and denote by $\{B_t\}_{t\geq 0}$ the standard Brownian motion starting from $0$. We have the following theorem.

\begin{theorem}\label{theorem 2.1 CLT}
Assuming that $d\geq 2$ and $\eta_0$ is distributed with $\nu_p$. Let $F: \{0, 1\}^{\mathcal{T}_d}\rightarrow \mathbb{R}$ be a local function satisfying \eqref{equ basic mean zero assumption} and $\xi_t^F$ be defined as in \eqref{equ 1.2 additive functional}, then
\[
\left\{\frac{1}{\sqrt{N}}\xi_{tN}^F:~0\leq t\leq T\right\}
\]
converges weakly, with respect to the uniform topology of $C[0, T]$, to $\{\sigma_FB_t\}_{0\leq t\leq T}$ as $N\rightarrow+\infty$, where
$\sigma_F$ is defined as in Lemma \ref{lemma 1.1}.
\end{theorem}

\begin{remark}\label{remark 2.1}
As we have introduced in Section \ref{section one}, Reference \cite{Sethuraman1996} proves the CLT of the additive functional of the SSEP in the case of $\mathbb{Z}^d$ by giving a sharp estimation of the spectral gap of the process. In this paper, to deal with the regular tree case, we utilize a different method. By applying the graphical representation of the SSEP introduced in \cite{Faggionato2024}, we extend the martingale decomposition strategy introduced in \cite{Kipnis1987} to the case of general additive functionals. For mathematical details, see Section \ref{section four}.
\end{remark}

\begin{remark}\label{remark 2.2}
By Theorem \ref{theorem 2.1 CLT}, the scaling function of the CLT is $\sqrt{N}$ for all $d\geq 2$, i.e., there is no $d$-dependent phase transition. This conclusion is natural according to Lemma \ref{lemma 1.1}, which shows that ${\rm Var}(\xi_{tN}^F)=O(N)$. Lemma \ref{lemma 1.1} relies on the fact that the simple random walk on $\mathcal{T}_d$ is transient for all $d\geq 2$. For mathematical details, see Section \ref{section three}.
\end{remark}

Our second main result is about the MDP of the additive functional of the SSEP on $\mathcal{T}_d$.
\begin{theorem}\label{theorem 2.2 mdp}
Let $\{a_t\}_{t\geq 0}$ be a positive sequence such that
\[
\lim_{t\rightarrow+\infty}\frac{a_t}{t}=\lim_{t\rightarrow+\infty}\frac{\sqrt{t}}{a_t}=0.
\]
For any local function $F$ satisfying \eqref{equ basic mean zero assumption} and any closed set $C\subseteq \mathbb{R}$,
\begin{equation*}
\limsup_{t\rightarrow+\infty}\frac{t}{a_t^2}\log \mathbb{P}_{\nu_p}\left(\frac{1}{a_t}\xi_t^F\in C\right)
\leq -\inf_{u\in C}\frac{u^2}{2\sigma_F^2}.
\end{equation*}
For any open set $O\subseteq \mathbb{R}$,
\begin{equation*}
\liminf_{t\rightarrow+\infty}\frac{t}{a_t^2}\log \mathbb{P}_{\nu_p}\left(\frac{1}{a_t}\xi_t^F\in O\right)\geq -\inf_{u\in O}\frac{u^2}{2\sigma_F^2}.
\end{equation*}
\end{theorem}

We further investigate the sample path MDP of the additive functional. To state our result, we denote by $\hat{C}[0, T]$ the set of $f\in C[0, T]$ such that $f_0=0$. For any $f\in \hat{C}[0, T]$, we define
\[
I(f)=
\begin{cases}
\frac{1}{2\sigma_F^2}\int_0^T \left(f^\prime_t\right)^2dt & \text{~if~}f\text{~is absolutely continuous,}\\
+\infty & \text{~else}.
\end{cases}
\]
Note that, $I(\cdot)$ has the equivalent definition
\begin{equation}\label{equ equivalent definition of I}
I(f)=\sup_{g\in C^1[0, T]}\left\{f_Tg_T-\int_0^Tf_sg_s^\prime ds-\frac{1}{2}\int_0^T \sigma_F^2g_s^2ds\right\}
\end{equation}
for any $f\in \hat{C}[0, T]$. The equivalence of above two definitions of $I(\cdot)$ follows from a routine analysis utilizing the Riesz representation theorem, which we omit here.

We have the following theorem.
\begin{theorem}\label{theorem 2.3 sample path mdp}
Let $\{b_N\}_{N\geq 1}$ be a positive sequence such that
\[
\lim_{N\rightarrow+\infty}\frac{b_N}{N}=\lim_{N\rightarrow+\infty}\frac{\sqrt{N\log N}}{b_N}=0.
\]
For any local function $F$ satisfying \eqref{equ basic mean zero assumption} and any closed set $\mathcal{C}\subseteq \hat{C}[0, T]$,
\begin{equation}\label{equ sample path MDP upper bound}
\limsup_{N\rightarrow+\infty}\frac{N}{b_N^2}\log \mathbb{P}_{\nu_p}\left(\left\{\frac{1}{b_N}\xi_{tN}^F:~0\leq t\leq T\right\}\in \mathcal{C}\right)
\leq -\inf_{f\in \mathcal{C}}I(f).
\end{equation}
For any open set $\mathcal{O}\subseteq \hat{C}[0, T]$,
\begin{equation}\label{equ sample path MDP lower bound}
\liminf_{N\rightarrow+\infty}\frac{N}{b_N^2}\log \mathbb{P}_{\nu_p}\left(\left\{\frac{1}{b_N}\xi_{tN}^F:~0\leq t\leq T\right\}\in \mathcal{O}\right)
\geq -\inf_{f\in \mathcal{O}}I(f).
\end{equation}
\end{theorem}
We believe that Theorem \ref{theorem 2.3 sample path mdp} holds for any $b_N$ between $\sqrt{N}$ and $N$. However, the strategy of our proof requires that $N=o(e^{b_N^2/N})$ and hence currently we need $b_N\gg \sqrt{N\log N}$. By Cauchy-Schwartz inequality, it is easy to check that
\[
\inf\left\{\int_0^1\left(f_s^\prime\right)^2ds:~f\in \hat{C}[0, 1], f(1)=u\right\}=u^2.
\]
Hence, for any $a_t$ between $\sqrt{t\log t}$ and $t$, Theorem \ref{theorem 2.2 mdp} is a corollary of Theorem \ref{theorem 2.3 sample path mdp} according to the contraction principle. Note that we still list Theorem \ref{theorem 2.2 mdp} as a main result since the case where $a_t$ between $\sqrt{t}$ and $\sqrt{t\log t}$ is not a direct consequence of Theorem \ref{theorem 2.3 sample path mdp}.

The proofs of Theorems \ref{theorem 2.2 mdp} and \ref{theorem 2.3 sample path mdp} utilize the exponential martingale strategy introduced in \cite{Kipnis1989}, where our extended version of the martingale decomposition formula introduced in \cite{Kipnis1987} still plays the key role as in the proof of the CLT.

The remainder of this paper is arranged as follows. In Section \ref{section three}, we recall the graphical representation of the SSEP introduced in \cite{Faggionato2024}. By utilizing this graphical representation, we prove Lemma \ref{lemma 1.1} and extend the martingale decomposition formula given in \cite{Kipnis1987}. According to our extended formula, we decompose the additive functional as a martingale plus a remainder. In Section \ref{section four}, we prove Theorem \ref{theorem 2.1 CLT}, i.e., the CLT. We show that the remainder term in the above decomposition converges weakly to $0$ via a calculation of the second moment and then show that the martingale term converges weakly to the target Brownian motion via a calculation of the quadratic variation process. In Sections \ref{section five} and \ref{section six}, we prove Theorems \ref{theorem 2.2 mdp} and \ref{theorem 2.3 sample path mdp}, i.e., the MDP and the sample path MDP. As we have introduced aboved, the proofs of these two MDPs utilize the exponential martingale strategy introduced in \cite{Kipnis1989} and our extended martingale decomposition formula, where two replacement lemmas play the key role. To prove these replacement lemmas, we utilize a variation formula for the largest eigenvalue of self-adjoint operator given in \cite{kipnis+landim99} and a time version block-analysis introduced in \cite{Xue2025}.

\section{Graphical representation}\label{section three}
In this section, we recall the graphical representation of the SSEP introduced in \cite{Faggionato2024}. As applications of this graphical representation, we prove Lemma \ref{lemma 1.1} and give an extended version of the martingale decomposition formula given in \cite{Kipnis1987}. According to the second explanation of the definition of $\{\eta_t\}_{t\geq 0}$ given in Section \ref{section one}, for each pair of neighbors $x, y\in \mathcal{T}_d$, the event moments at which $x$ and $y$ exchange their values with each other form a Poisson flow with rate $1$. We denote this Poisson flow by $\{N_t^{x, y}\}_{t\geq 0}$. Consequently, when $\eta_0$ is given, $\eta_t$ is determined by $\eta_0$ and $\{N_s^{x,y}:~0\leq s\leq t, x\sim y\}$. In detail, for any $0\leq s\leq t$ and $x\in \mathcal{T}_d$, there exists an unique $X_s^{t, x}\in \mathcal{T}_d$ such that $\eta_t(x)=\eta_{t-s}(X_s^{t, x})$ and especially,
\begin{equation}\label{equ 3.1 graphical representation}
\eta_t(x)=\eta_0(X_t^{t, x}).
\end{equation}
The random vertex $X_\cdot^{t, x}$ jumps according to the following rule. At moment $0$, $X_0^{t, x}=x$. For any $0\leq s_1<s_2\leq t$, we denote by $z$ the vertex $X_{s_1}^{t, x}$. If there is no event moment of $N_\cdot^{z, y}$ in $[t-s_2, t-s_1]$ for any $y\sim z$, then $X_r^{t, x}=z$ for all $r\in [s_1, s_2]$. For any $0\leq s\leq t$, let $X_{s-}^{t, x}$ be the position of $X_\cdot^{t, x}$ at the moment just before $s$. If $t-s$ is an event moment of $N_{\cdot}^{X_{s-}^{t, x}, y}$ for some $y\sim X_{s-}^{t, x}$, then $X_s^{t, x}=y$. In conclusion, $\{X_s^{t, x}\}_{0\leq s\leq t}$ is a simple random walk on $\mathcal{T}_d$ with generator $\mathcal{G}$ given by
\[
\mathcal{G}f(x)=\sum_{y\sim x}\left(f(y)-f(x)\right)
\]
for any $x\in \mathcal{T}_d$ and $f: \mathcal{T}_d\rightarrow \mathbb{R}$.

The coupling relationship \eqref{equ 3.1 graphical representation} is called a `graphical representation' since this formula can also be explained by a graphical method. In detail, we consider the graph $\mathcal{T}_d\times [0, +\infty)$, i.e., there is a time axis at each $x\in \mathcal{T}_d$. At each event moment $s$ of $N_\cdot^{x, y}$, we write a double-headed arrow `$\leftrightarrow$' connecting $(x, s)$ and $(y, s)$. Since all vertices constantly exchange their values with neighbors, we can backtrack where the value of $x$ at moment $t$ comes from. For any $0\leq r\leq t$, let $X_r^{t, x}$ be the position which provides the value $\eta_t(x)$ at moment $t-r$. Hence, if there is an $\leftrightarrow$ connecting $(X_{s-}^{t, x}, t-s)$ and $(y, t-s)$ for some $y\sim X_{s-}^{t, x}$, then $X_\cdot^{t, x}$ jumps from $X_{s-}^{t, x}$ to $y$ at moment $s$ through this $\leftrightarrow$, i.e., the provider $X_{s-}^{t, x}$ of $\eta_t(x)$ is replaced by $y$.

For integer $m\geq 1$ and different $x_1, x_2, \ldots, x_m\in \mathcal{T}_d$, we discuss the joint distribution of $\left\{\left(X_s^{t, x_1},\ldots, X_s^{t, x_m}\right)\right\}_{0\leq s\leq t}$. For $1\leq i\neq j\leq m$, if $X_{s-}^{t, x_i}\sim X_{s-}^{t, x_j}$ and there is an $\leftrightarrow$ connecting
$\left(X_{s-}^{t, x_i}, t-s\right)$ and $\left(X_{s-}^{t, x_j}, t-s\right)$, then
\[
\left(X_{s}^{t, x_i}, X_{s}^{t, x_j}\right)=\left(X_{s-}^{t, x_j}, X_{s-}^{t, x_i}\right)
\]
according to the evolution of $X_\cdot^{t, x}$. Hence, $\left\{\left(X_s^{t, x_1},\ldots, X_s^{t, x_m}\right)\right\}_{0\leq s\leq t}$ performs independent simple random walks on $\mathcal{T}_d$ except that when any two components are neighbors, these two components exchange their positions with each other at rate $1$ to avoid collision. As a result, for each integer $m\geq 1$, $x_1, \ldots, x_m, y_1,\ldots, y_m\in \mathcal{T}_d$ and any $t>s$, $\mathbb{P}\left(\left(X_s^{t, x_1},\ldots, X_s^{t, x_m}\right)=(y_1,\ldots, y_m)\right)$ is independent of the choice of $t$. Hence, we denote $\mathbb{P}\left(\left(X_s^{t, x_1},\ldots, X_s^{t, x_m}\right)=(y_1,\ldots, y_m)\right)$ by $Q_s^m\left((x_1, \ldots, x_m), (y_1, \ldots, y_m)\right)$.
Especially, when $m=1$, $\{Q_s^1(\cdot, \cdot)\}_{0\leq s\leq t}$ are the transition probabilities of the simple random walk on $\mathcal{T}_d$ as we have explained above.

Now we prove Lemma \ref{lemma 1.1}.

\proof[Proof of Lemma \ref{lemma 1.1}]

Since $F$ is a local function, there exist integer $m\geq 1$, $x_1, \ldots, x_m\in \mathcal{T}_d$ which are different with each other and $H:\{0, 1\}^m\rightarrow\mathbb{R}$ such that
\[
F(\eta)=H\left(\eta(x_1),\ldots, \eta(x_m)\right).
\]
Then, by \eqref{equ 3.1 graphical representation}, we have
\begin{equation}\label{equ 3.2 mean calculation}
\mathbb{E}_\eta F(\eta_t)=\sum_{y_1,\ldots, y_m}Q_t^m\left((x_1,\ldots,x_m), (y_1,\ldots,y_m)\right)H\left(\eta(y_1),\ldots, \eta(y_m)\right)
\end{equation}
for any $t\geq 0$, where the sum is over all $(y_1, \ldots, y_m)$ such that $y_i\neq y_j$ for all $1\leq i<j\leq m$. From now on, for simplicity, we write $Q_t^m\left((x_1,\ldots,x_m), (y_1,\ldots,y_m)\right)$ as $Q_{t, x_1,\ldots,x_m}^m(y_1,\ldots, y_m)$ and $H\left(\eta(y_1),\ldots, \eta(y_m)\right)$ as $H(\eta, y_1,\ldots, y_m)$.

According to the bilinear property of the covariance operator and the invariance of $\nu_p$, we have
\begin{align*}
{\rm Var}_{\nu_p}\left(\xi_t^F\right)&=2\int_0^t\left(\int_0^s {\rm Cov}_{\nu_p}\left(F(\eta_u), F(\eta_0)\right)du\right)ds\\
&=2\int_0^t\left(\int_0^s \mathbb{E}_{\nu_p}\left(F(\eta_u)F(\eta_0)\right)du\right)ds.
\end{align*}
Note that the last equation in the above formula follows from Assumption \eqref{equ basic mean zero assumption}. Hence,
\[
\lim_{t\rightarrow+\infty}\frac{1}{t}{\rm Var}_{\nu_p}\left(\xi_t^F\right)=\sigma_F^2,
\]
where
\[
\sigma_F^2=2\int_0^{+\infty}\mathbb{E}_{\nu_p}\left(F(\eta_u)F(\eta_0)\right)du.
\]
We only need to show that $\sigma_F^2<+\infty$. According to \eqref{equ 3.2 mean calculation} and the Markov property of $\{\eta_t\}_{t\geq 0}$, we have
\begin{align}\label{equ 3.2 two}
&\mathbb{E}_{\nu_p}\left(F(\eta_u)F(\eta_0)\right)\\
&=\sum_{y_1,\ldots,y_m}\mathbb{E}_{\nu_p}\left(H(\eta_0, x_1,\ldots, x_m)H(\eta_0, y_1,\ldots, y_m)\right)Q_{u, x_1, \ldots, x_m}^m(y_1,\ldots, y_m).\notag
\end{align}
We denote by $\Lambda_{x_1,\ldots, x_m}$ the set of all $(y_1, \ldots, y_m)$ such that $\{y_1,\ldots, y_m\}\bigcap \{x_1, \ldots, x_m\}\neq \emptyset$. According to the definition of $\nu_p$ and Assumption \eqref{equ basic mean zero assumption},
\[
\mathbb{E}_{\nu_p}\left(H(\eta_0, x_1,\ldots, x_m)H(\eta_0, y_1,\ldots, y_m)\right)=0
\]
for $(y_1,\ldots, y_m)\not\in \Lambda_{x_1,\ldots, x_m}$. Let $K_H=\max_{\vec{w}\in \{0, 1\}^m}\left|H(\vec{w})\right|$, then
\begin{align}\label{equ 3.3}
\mathbb{E}_{\nu_p}\left(F(\eta_u)F(\eta_0)\right)&\leq K_H^2\sum_{(y_1,\ldots,y_m)\in \Lambda_{x_1,\ldots,x_m}}Q_{u, x_1, \ldots, x_m}^m(y_1,\ldots, y_m) \notag\\
&= K_H^2\mathbb{P}\left(\left\{X_u^{t, x_1},\ldots, X_u^{t, x_m}\right\}\bigcap\left\{x_1,\ldots, x_m\right\}\neq \emptyset\right)\notag\\
&\leq K_H^2\sum_{i=1}^m\sum_{j=1}^mQ^1_{u, x_i}(x_j).
\end{align}
It is shown in \cite{Xue2025} that
\begin{equation}\label{equ random walk heat estimation}
Q^1_{u, x}(z)\leq e^{-u(\sqrt{d}-1)^2}
\end{equation}
for any $x, z\in \mathcal{T}_d$. Consequently, $\sigma_F^2<+\infty$ follows directly from \eqref{equ 3.3}.
\qed

Now we extend the martingale decomposition formula introduced in \cite{Kipnis1987} to the case of general additive functionals. For a local function $F$ satisfying \eqref{equ basic mean zero assumption}, as in the proof of Lemma \ref{lemma 1.1}, we can write $F(\eta)$ as $H(\eta, x_1, \ldots, x_m)$ for some integer $m\geq 1$, $x_1,\ldots, x_m\in \mathcal{T}_d$ and $H: \{0, 1\}^m\rightarrow\mathbb{R}$. Note that $H(\eta, x_1, \ldots, x_m)=H\left(\eta(x_1),\ldots, \eta(x_m)\right)$ defined as in the proof of Lemma \ref{lemma 1.1}. For any $\lambda>0$ and $y_1,\ldots, y_m\in \mathcal{T}_d$, we define
\[
\beta_{\lambda, x_1,\ldots,x_m}(y_1,\ldots, y_m)=\int_0^{+\infty}e^{-\lambda s}Q^m_{s, x_1,\ldots, x_m}(y_1,\ldots, y_m)ds.
\]
Then, we define
\[
G_\lambda^F(\eta)=\sum_{y_1,\ldots, y_m}H(\eta, y_1, \ldots, y_m)\beta_{\lambda, x_1,\ldots,x_m}(y_1,\ldots, y_m).
\]
According to the definition of $\mathcal{L}$,
\begin{align*}
&\mathcal{L}G_\lambda^F(\eta)\\
&=\frac{1}{2}\sum_{y\in \mathcal{T}_d}\sum_{z\sim y}\left(\sum_{y_1,\ldots, y_m}\left(H(\eta^{y, z}, y_1, \ldots, y_m)-H(\eta, y_1,\ldots, y_m)\right)\right)\beta_{\lambda, x_1,\ldots, x_m}(y_1,\ldots, y_m).
\end{align*}
We denote by $\gamma_{y_1,\ldots, y_m}$ the set of all $(z_1, \ldots, z_m)$ such that $\left(X_\cdot^{t, x_1},\ldots, X_\cdot^{t, x_m}\right)$ jumps from $(y_1,\ldots, y_m)$ to $(z_1, \ldots, z_m)$ at rate $1$. If $y\sim z$ make $H(\eta^{y, z}, y_1, \ldots, y_m)-H(\eta, y_1,\ldots, y_m)\neq 0$, then
$\{y, z\}\bigcap \{y_1,\ldots, y_m\}\neq \emptyset$ and
\[
H(\eta^{y, z}, y_1, \ldots, y_m)=H(\eta, z_1, \ldots, z_m)
\]
for some $(z_1, \ldots, z_m)\in \gamma_{y_1, \ldots, y_m}$. According to the Kolmogorov-Chapman equation,
\begin{align*}
&\frac{d}{ds}Q^m_{s, x_1,\ldots, x_m}(y_1,\ldots, y_m)\\
&=\sum_{(z_1,\ldots, z_m)\in \gamma_{y_1, \ldots, y_m}}\left(Q^m_{s, x_1,\ldots, x_m}(z_1,\ldots, z_m)-Q^m_{s, x_1,\ldots, x_m}(y_1,\ldots, y_m)\right)
\end{align*}
and hence, by the integration-by-parts formula,
\begin{align}\label{equ 3.4 generator on G}
&\mathcal{L}G_\lambda^F(\eta)\notag\\
&=\sum_{y_1,\ldots, y_m}\sum_{(z_1, \ldots, z_m)\in \gamma_{y_1, \ldots, y_m}}\left(H(\eta, z_1, \ldots, z_m)-H(\eta, y_1, \ldots, y_m)\right)\beta_{\lambda, x_1,\ldots, x_m}(y_1,\ldots, y_m)\notag\\
&=\sum_{y_1,\ldots, y_m}H(\eta, y_1, \ldots, y_m)\left(\sum_{z_1, \ldots, z_m}\left(\beta_{\lambda, x_1,\ldots, x_m}(z_1,\ldots, z_m)-\beta_{\lambda, x_1,\ldots, x_m}(y_1,\ldots, y_m)\right)\right)\notag\\
&=\sum_{y_1, \ldots, y_m}H(\eta, y_1, \ldots, y_m)\int_0^{+\infty}e^{-\lambda s}\frac{d}{ds}Q^m_{s, x_1,\ldots, x_m}(y_1,\ldots, y_m)ds\notag\\
&=\sum_{y_1, \ldots, y_m}H(\eta, y_1, \ldots, y_m)\left(-Q^m_{0, x_1,\ldots, x_m}(y_1,\ldots, y_m)+\lambda\int_0^{+\infty}e^{-\lambda s}Q^m_{s, x_1,\ldots, x_m}(y_1,\ldots, y_m)ds\right)\notag\\
&=-H(\eta, x_1,\ldots, x_m)+\lambda G_\lambda^F(\eta)=\lambda G_\lambda^F(\eta)-F(\eta).
\end{align}
We define
\[
M_t^{F, \lambda}=G_\lambda^F(\eta_t)-G_\lambda^F(\eta_0)-\int_0^t \mathcal{L}G_\lambda^F(\eta_s)ds,
\]
then $\{M_t^{F, \lambda}\}_{t\geq 0}$ is a martingale according to the Dynkin's martingale formula. By \eqref{equ 3.4 generator on G}, we have
\begin{equation}\label{equ 3.5 extended martingale decomposition formula}
\xi_t^F=M_t^{F, \lambda}-\left(G_\lambda^F(\eta_t)-G_\lambda^F(\eta_0)-\int_0^t \lambda G_\lambda^F(\eta_s)ds\right).
\end{equation}
When $F(\eta)=\eta(x)-p$, Equation \eqref{equ 3.5 extended martingale decomposition formula} reduces to the martingale decomposition formula introduced in \cite{Kipnis1987}, which is utilized in the proof of the occupation time CLT.

\section{Proof of Theorem \ref{theorem 2.1 CLT}}\label{section four}
In this section, we prove Theorem \ref{theorem 2.1 CLT}. By \eqref{equ 3.5 extended martingale decomposition formula},
\[
\frac{1}{\sqrt{N}}\xi_{tN}^F=\frac{1}{\sqrt{N}}M_{tN}^{F, 1/N}-\frac{1}{\sqrt{N}}\left(G_{1/N}^F(\eta_{tN})-G_{1/N}^F(\eta_0)-\frac{1}{N}\int_0^{tN} G_{1/N}^F(\eta_s)ds\right).
\]
The following lemma shows that the remainder term
\[
\frac{1}{\sqrt{N}}\left(G_{1/N}^F(\eta_{tN})-G_{1/N}^F(\eta_0)-\frac{1}{N}\int_0^{tN} G_{1/N}^F(\eta_s)ds\right)
\]
in the above martingale decomposition formula converges weakly to $0$.

\begin{lemma}\label{lemma 4.1}
For any $t\geq 0$,
\[
\lim_{N\rightarrow+\infty}\frac{1}{\sqrt{N}}\left(G_{1/N}^F(\eta_{tN})-G_{1/N}^F(\eta_0)-\frac{1}{N}\int_0^{tN} G_{1/N}^F(\eta_s)ds\right)=0
\]
in $L^2$.
\end{lemma}

\proof

According to the invariance of $\nu_p$ and Cauchy-Schwartz inequality, we only need to show that
\begin{equation}\label{equ 4.1}
\lim_{N\rightarrow+\infty}\mathbb{E}_{\nu_p}\left(\left(\frac{1}{\sqrt{N}}G_{1/N}^F(\eta_0)\right)^2\right)=0.
\end{equation}
According to the definition of $G_\lambda^F$, we have
\begin{align*}
\mathbb{E}_{\nu_p}\left(\left(\frac{1}{\sqrt{N}}G_{1/N}^F(\eta_0)\right)^2\right)=\frac{1}{N}\sum_{y_1,\ldots,y_m}\sum_{z_1,\ldots, z_m}\mathcal{A}_{x_1, \ldots, x_m}(N, y_1, \ldots, y_m, z_1,\ldots, z_m),
\end{align*}
where
\begin{align*}
&\mathcal{A}_{x_1, \ldots, x_m}(N, y_1, \ldots, y_m, z_1,\ldots, z_m)\\
&=\mathbb{E}_{\nu_p}\left(H(\eta, y_1,\ldots, y_m)H(\eta, z_1, \ldots, z_m)\right)
\beta_{1/N, x_1,\ldots,x_m}(y_1,\ldots, y_m)\beta_{1/N, x_1,\ldots,x_m}(z_1,\ldots, z_m).
\end{align*}
Let $\Lambda_{y_1,\ldots, y_m}$ be the set of all $(z_1,\ldots, z_m)$ such that $\{z_1,\ldots, z_m\}\bigcap \{y_1,\ldots, y_m\}\neq\emptyset$. For any $(v_1,\ldots, v_m)\not\in \Lambda_{y_1,\ldots, y_m}$,
\[
\mathcal{A}_{x_1, \ldots, x_m}(N, y_1, \ldots, y_m, v_1,\ldots, v_m)=0
\]
according to the definition of $\nu_p$ and Assumption \eqref{equ basic mean zero assumption}. Hence,
\begin{align*}
&\mathbb{E}_{\nu_p}\left(\left(\frac{1}{\sqrt{N}}G_{1/N}^F(\eta_0)\right)^2\right)\\
&=\frac{1}{N}\sum_{y_1,\ldots,y_m}\sum_{(z_1,\ldots, z_m)\in \Lambda_{y_1,\ldots, y_m}}\mathcal{A}_{x_1, \ldots, x_m}(N, y_1, \ldots, y_m, z_1,\ldots, z_m)\\
&\leq \frac{K_H^2}{N}\sum_{y_1,\ldots,y_m}\sum_{(z_1,\ldots, z_m)\in \Lambda_{y_1,\ldots, y_m}}\beta_{1/N, x_1,\ldots,x_m}(y_1,\ldots, y_m)\beta_{1/N, x_1,\ldots,x_m}(z_1,\ldots, z_m)\\
&\leq \frac{K_H^2}{N}\sum_{y_1,\ldots,y_m}\sum_{(z_1,\ldots, z_m)\in \Lambda_{y_1,\ldots, y_m}}\int_0^{+\infty}\int_0^{+\infty}Q^m_{s, x_1,\ldots, x_m}(y_1,\ldots, y_m)Q^m_{u, x_1,\ldots, x_m}(z_1,\ldots, z_m)dsdu,
\end{align*}
where $K_H=\max_{\vec{w}\in \{0, 1\}^m}|H(\vec{w})|$ defined as in Section \ref{section three}. Let $\left\{\left(\hat{X}_r^{t, x_1},\ldots, \hat{X}_r^{t, x_m}\right)\right\}_{0\leq r\leq t}$ be an independent copy of $\left\{\left(X_r^{t, x_1},\ldots, X_r^{t, x_m}\right)\right\}_{0\leq r\leq t}$, then for $u,s\leq t$,
\begin{align*}
&\sum_{y_1,\ldots,y_m}\sum_{(z_1,\ldots, z_m)\in \Lambda_{y_1,\ldots, y_m}}Q^m_{s, x_1,\ldots, x_m}(y_1,\ldots, y_m)Q^m_{u, x_1,\ldots, x_m}(z_1,\ldots, z_m)\\
&=\mathbb{P}\left(\left\{\hat{X}_s^{t, x_1},\ldots, \hat{X}_s^{t, x_m}\right\}\bigcap \left\{X_u^{t, x_1},\ldots, X_u^{t, x_m}\right\}\neq \emptyset\right)\\
&\leq \sum_{i=1}^m\sum_{j=1}^m\mathbb{P}\left(\hat{X}_s^{t, x_i}=X_u^{t, x_j}\right)=\sum_{i=1}^m\sum_{j=1}^m\sum_{x\in \mathcal{T}_d}Q_{s, x_i}^1(x)Q_{u, x_j}^1(x)\\
&=\sum_{i=1}^m\sum_{j=1}^m Q_{s+u, x_i}^1(x_j).
\end{align*}
Since we can choose $t$ arbitrarily large, the above inequality holds for all $s,u\geq 0$. In conclusion,
\[
\mathbb{E}_{\nu_p}\left(\left(\frac{1}{\sqrt{N}}G_{1/N}^F(\eta_0)\right)^2\right)
\leq \frac{K_H^2}{N}\sum_{i=1}^m\sum_{j=1}^m \int_0^{+\infty}\int_0^{+\infty}Q_{s+u, x_i}^1(x_j)dsdu.
\]
Then, by \eqref{equ random walk heat estimation}, $\mathbb{E}_{\nu_p}\left(\left(\frac{1}{\sqrt{N}}G_{1/N}^F(\eta_0)\right)^2\right)=O(N^{-1})$ and the proof is complete.
\qed

The next lemma shows that the martingale part $\frac{1}{\sqrt{N}}M_{tN}^{F, 1/N}$ in our decomposition formula converges weakly to the target Brownian motion.

\begin{lemma}\label{lemma 4.2 martingale part convergence}
As $N\rightarrow+\infty$,
\[
\left\{\frac{1}{\sqrt{N}}M_{tN}^{F, 1/N}:~0\leq t\leq T\right\}
\]
converges weakly, with respect to the Skorohod topology of $D[0, T]$, to $\{\sigma_FB_t\}_{0\leq t\leq T}$.
\end{lemma}

\proof

For any $\lambda>0$, let
\[
J_t^{F, \lambda}=\int_0^t\mathcal{L}\left(\left(G_\lambda^F(\eta_s)\right)^2\right)-2G_\lambda^F(\eta_s)\mathcal{L}G_\lambda^F(\eta_s)ds.
\]
According to the Dynkin's martingale formula, $\left\{\left(M_t^{F, \lambda}\right)^2-J_t^{F, \lambda}\right\}_{t\geq 0}$ is a martingale. By Theorem 1.4 in Chapter 7 of \cite{Ethier1986}, to complete the proof we only need to show that
\begin{equation}\label{equ 4.2 quadratic variation process convergence}
\lim_{N\rightarrow+\infty}\frac{1}{N}J_{tN}^{F, 1/N}=\sigma_F^2t
\end{equation}
in $L^2$ for all $t\geq 0$. To prove \eqref{equ 4.2 quadratic variation process convergence}, we only need to show that
\begin{equation}\label{equ 4.3}
\lim_{N\rightarrow+\infty}\mathbb{E}_{\nu_p}\left(\frac{1}{N}J_{tN}^{F, 1/N}\right)=\sigma_F^2t
\end{equation}
and
\begin{equation}\label{equ 4.4}
\lim_{N\rightarrow+\infty}{\rm Var}_{\nu_p}\left(\frac{1}{N}J_{tN}^{F, 1/N}\right)=0.
\end{equation}
We first check \eqref{equ 4.3}. According to the definition of $\mathcal{L}$,
\[
\left(\left(G_{1/N}^F(\eta)\right)^2\right)-2G_{1/N}^F(\eta)\mathcal{L}G_{1/N}^F(\eta)
=\frac{1}{2}\sum_{y\in \mathcal{T}_d}\sum_{z\sim y}\left(G_{1/N}^F(\eta^{y, z})-G_{1/N}^F(\eta)\right)^2.
\]
Hence, according to the invariance of $\nu_p$,
\[
\mathbb{E}_{\nu_p}\left(\frac{1}{N}J_{tN}^{F, 1/N}\right)=\frac{t}{2}\mathbb{E}_{\nu_p}\left(\sum_{y\in \mathcal{T}_d}\sum_{z\sim y}\left(G_{1/N}^F(\eta^{y, z})-G_{1/N}^F(\eta)\right)^2\right).
\]
According to the fact that $\eta\rightarrow \eta^{y, z}$ is a measure preserving transformation under $\nu_p$, we have
\begin{align*}
&\mathbb{E}_{\nu_p}\left(\sum_{y\in \mathcal{T}_d}\sum_{z\sim y}\left(G_{1/N}^F(\eta^{y, z})-G_{1/N}^F(\eta)\right)^2\right)\\
&=2\mathbb{E}_{\nu_p}\left(\sum_{y_1, \ldots, y_m}\sum_{z_1,\ldots, z_m}\mathcal{B}(\eta, N, y_1, \ldots, y_m, z_1,\ldots, z_m)\right),
\end{align*}
where
\begin{align*}
&\mathcal{B}(\eta, N, y_1, \ldots, y_m, z_1,\ldots, z_m)\\
&=H(z_1,\ldots, z_m)\beta_{1/N, x_1,\ldots, x_m}(z_1,\ldots, z_m)\\
&\text{\quad\quad}\times\sum_{y\in \mathcal{T}_d}\sum_{z\sim y}\left(H(\eta, y_1,\ldots, y_m)-H(\eta^{y, z}, y_1, \ldots, y_m)\right)\beta_{1/N, x_1,\ldots, x_m}(y_1,\ldots, y_m).
\end{align*}
According to an analysis similar with that leading to \eqref{equ 3.5 extended martingale decomposition formula},
\begin{align*}
&\sum_{y_1,\ldots, y_m}\sum_{y\in \mathcal{T}_d}\sum_{z\sim y}\left(H(\eta, y_1,\ldots, y_m)-H(\eta^{y, z}, y_1, \ldots, y_m)\right)\beta_{1/N, x_1,\ldots, x_m}(y_1,\ldots, y_m)\\
&=2\sum_{y_1, \ldots, y_m}H(\eta, y_1, \ldots, y_m)\\
&\text{\quad\quad\quad}\times\left(\sum_{(v_1,\ldots, v_m)\in \gamma_{y_1, \ldots, y_m}}\left(\beta_{1/N, x_1,\ldots, x_m}(y_1,\ldots, y_m)-\beta_{1/N, x_1,\ldots, x_m}(v_1,\ldots, v_m)\right)\right)\\
&=-2\sum_{y_1,\ldots, y_m}H(\eta, y_1,\ldots, y_m)\int_0^{+\infty}e^{-s/N}\frac{d}{ds}Q^m_{s, x_1,\ldots, x_m}(y_1,\ldots, y_m)ds\\
&=2H(\eta, x_1, \ldots, x_m)-\frac{2}{N}\sum_{y_1,\ldots, y_m}H(\eta, y_1,\ldots, y_m)\int_0^{+\infty}e^{-s/N}Q^m_{s, x_1,\ldots, x_m}(y_1,\ldots, y_m)ds.
\end{align*}
Hence,
\[
\mathbb{E}_{\nu_p}\left(\sum_{y_1, \ldots, y_m}\sum_{z_1,\ldots, z_m}\mathcal{B}(\eta, N, y_1, \ldots, y_m, z_1,\ldots, z_m)\right)
={\rm \uppercase\expandafter{\romannumeral1}}+{\rm \uppercase\expandafter{\romannumeral2}},
\]
where
\[
{\rm \uppercase\expandafter{\romannumeral1}}=2\mathbb{E}_{\nu_p}\left(\sum_{z_1,\ldots, z_m}H(\eta, z_1,\ldots, z_m)H(\eta, x_1,\ldots, x_m)\beta_{1/N, x_1,\ldots, x_m}(z_1,\ldots, z_m)\right)
\]
and
\begin{align*}
{\rm \uppercase\expandafter{\romannumeral2}}&=-\frac{2}{N}\sum_{z_1,\ldots, z_m}\sum_{y_1,\ldots, y_m}\mathbb{E}_{\nu_p}\left(H(\eta, z_1, \ldots, z_m)H(\eta, y_1, \ldots, y_m)\right)\\
&\text{\quad\quad\quad\quad}\times\beta_{1/N, x_1,\ldots, x_m}(z_1,\ldots, z_m)\beta_{1/N, x_1,\ldots, x_m}(y_1,\ldots, y_m)\\
&=-2\mathbb{E}_{\nu_p}\left(\left(\frac{1}{\sqrt{N}}G_{1/N}^F(\eta_0)\right)^2\right).
\end{align*}
According to the definition of $\beta_{\lambda, x_1, \ldots, x_m}$ and \eqref{equ 3.2 two},
\begin{align}\label{equ 4.5}
\lim_{N\rightarrow+\infty}{\rm \uppercase\expandafter{\romannumeral1}}
&=2\mathbb{E}_{\nu_p}\left(\sum_{z_1,\ldots, z_m}H(\eta, z_1,\ldots, z_m)H(\eta, x_1,\ldots, x_m)\beta_{0, x_1,\ldots, x_m}(z_1,\ldots, z_m)\right) \notag\\
&=2\int_0^{+\infty}\sum_{z_1,\ldots, z_m}\mathbb{E}_{\nu_p}\left(H(\eta, z_1,\ldots, z_m)H(\eta, x_1,\ldots, x_m)\right)Q_{s, x_1, \ldots, x_m}^m(z_1,\ldots, z_m)ds\notag\\
&=2\int_0^{+\infty}\mathbb{E}_{\nu_p}\left(F(\eta_0)F(\eta_s)\right)ds=\sigma_F^2.
\end{align}
As we have shown in the proof of Lemma \ref{lemma 4.1},
\[
\lim_{N\rightarrow+\infty}{\rm \uppercase\expandafter{\romannumeral2}}
=-2\lim_{N\rightarrow+\infty}\mathbb{E}_{\nu_p}\left(\left(\frac{1}{\sqrt{N}}G_{1/N}^F(\eta_0)\right)^2\right)=0.
\]
Therefore, Equation \eqref{equ 4.3} follows from \eqref{equ 4.5}.

At last, we check \eqref{equ 4.4}. According to the definition of $G_\lambda^F$,
\begin{align*}
&G_{1/N}^F(\eta^{y, z})-G_{1/N}^F(\eta)\\
&=\sum_{(y_1,\ldots, y_n)\in \Lambda_{y, z}}\left(H(\eta^{y, z}, y_1,\ldots, y_m)-H(\eta, y_1,\ldots, y_m)\right)\beta_{1/N, x_1, \ldots, x_m}(y_1, \ldots, y_m),
\end{align*}
where $\Lambda_{y, z}$ is the set of $(y_1,\ldots, y_m)$ such that $\left\{y_1, \ldots, y_m\right\}\bigcap\{y, z\}\neq \emptyset$. We define
\[
\mathcal{K}(N, y, z)=\sum_{(y_1,\ldots, y_m)\in \Lambda_{y, z}}\sum_{(z_1, \ldots, z_m)\in \Lambda_{y, z}}\beta_{1/N, x_1,\ldots, x_m}(y_1,\ldots, y_m)\beta_{1/N, x_1, \ldots, x_m}(z_1, \ldots, z_m),
\]
then
\begin{align*}
\sum_{y\in \mathcal{T}_d}\sum_{z\sim y}\mathcal{K}(N, y, z)\leq \sum_{y\in \mathcal{T}_d}\sum_{z\sim y}\mathcal{K}(0, y, z)\leq \int_0^{+\infty}\int_0^{+\infty}\mathcal{R}(s, u)dsdu,
\end{align*}
where
\[
\mathcal{R}(s, u)=\mathbb{P}\left(D(X_s^{t, x_i}, \hat{X}_u^{t, x_j})\leq 1 \text{~for some~}1\leq i, j\leq m\right).
\]
Since each vertex on $\mathcal{T}_d$ has $d+1$ neighbors, according to \eqref{equ random walk heat estimation} and the total probability formula,
\[
\mathcal{R}(s, u)\leq m^2(d+2)e^{-\max\{s, u\}\left(\sqrt{d}-1\right)^2}.
\]
As a result,
\[
\sup_{N\geq 1}\left(\sum_{y\in \mathcal{T}_d}\sum_{z\sim y}\mathcal{K}(N, y, z)\right)^2<+\infty.
\]
Hence, according to the invariance of $\nu_p$ and the bilinear property of the covariance operator, to check \eqref{equ 4.4} we only need to show that
\begin{align}\label{equ 4.6}
\lim_{s\rightarrow+\infty}&\Bigg(\sup\Big|{\rm Cov}_{\nu_p}\big(\mathcal{U}(0, y, z, y_1, \ldots, y_m, z_1, \ldots, z_m),\notag\\
&\text{\quad\quad\quad\quad}\mathcal{U}(s, w, v, w_1, \ldots, w_m, v_1, \ldots, v_m)\big)\Big|\Bigg)=0,
\end{align}
where
\begin{align*}
&\mathcal{U}(s, w, v, w_1, \ldots, w_m, v_1, \ldots, v_m)\\
&=\left(H(\eta_s^{w, v}, w_1, \ldots, w_m)-H(\eta_s, w_1, \ldots, w_m)\right)\left(H(\eta_s^{w, v}, v_1, \ldots, v_m)-H(\eta_s, v_1, \ldots, v_m)\right)
\end{align*}
and the $\sup$ in \eqref{equ 4.6} is over all $y\sim z, w\sim v$, $(y_1,\ldots, y_m), (z_1, \ldots, z_m)\in \Lambda_{y, z}$ and
\[
(w_1, \ldots, w_m), (v_1, \ldots, v_m)\in \Lambda_{w, v}.
\]
By \eqref{equ 3.1 graphical representation}, the value of $H(\eta_s^{w, v}, w_1, \ldots, w_m)-H(\eta_s, w_1, \ldots, w_m)$ only depends on $\eta_0$ and
\[
X_s^{s, w_1}, \ldots, X_s^{s, w_m}, X_s^{s, w}, X_s^{s, v}.
\]
Hence, according to the definition of $\nu_p$, conditioned on
\[
\left(X_s^{s, w_1}, \ldots, X_s^{s, w_m}, X_s^{s, w}, X_s^{s, v}\right)=\left(\hat{w}_1, \ldots, \hat{w}_m, \hat{w}, \hat{v}\right)
\]
and
\[
\left(X_s^{s, v_1}, \ldots, X_s^{s, v_m}, X_s^{s, w}, X_s^{s, v}\right)=\left(\hat{v}_1, \ldots, \hat{v}_m, \hat{w}, \hat{v}\right)
\]
for some $\left(\hat{w}_1, \ldots, \hat{w}_m, \hat{w}, \hat{v}, \hat{v}_1, \ldots, \hat{v}_m\right)$ such that
\[
\left\{\hat{w}_1, \ldots, \hat{w}_m, \hat{w}, \hat{v}, \hat{v}_1, \ldots, \hat{v}_m\right\}\bigcap \{y_1, \ldots, y_m, z_1, \ldots, z_m, y, z\}=\emptyset,
\]
we have
\begin{align*}
{\rm Cov}_{\nu_p}\big(\mathcal{U}(0, y, z, y_1, \ldots, y_m, z_1, \ldots, z_m), \mathcal{U}(0, \hat{w}, \hat{v}, \hat{w}_1, \ldots, \hat{w}_m, \hat{v}_1, \ldots, \hat{v}_m)\big)=0.
\end{align*}
As a result,
\begin{align*}
&\left|{\rm Cov}_{\nu_p}\big(\mathcal{U}(0, y, z, y_1, \ldots, y_m, z_1, \ldots, z_m), \mathcal{U}(s, w, v, w_1, \ldots, w_m, v_1, \ldots, v_m)\big)\right|\\
&\leq 4K_H^2\left({\rm \uppercase\expandafter{\romannumeral3}}+{\rm \uppercase\expandafter{\romannumeral4}}\right),
\end{align*}
where
\[
{\rm \uppercase\expandafter{\romannumeral3}}=\mathbb{P}\left(\left\{X_s^{s, w_1}, \ldots, X_s^{s, w_m}, X_s^{s, w}, X_s^{s, v}\right\}
\bigcap\{y, z, y_1, \ldots, y_m, z_1, \ldots, z_m\}\neq \emptyset\right)
\]
and
\[
{\rm \uppercase\expandafter{\romannumeral4}}=\mathbb{P}\left(\left\{X_s^{s, v_1}, \ldots, X_s^{s, v_m}, X_s^{s, w}, X_s^{s, v}\right\}
\bigcap\{y, z, y_1, \ldots, y_m, z_1, \ldots, z_m\}\neq \emptyset\right).
\]
By \eqref{equ random walk heat estimation},
\begin{equation}\label{equ 4.7}
{\rm \uppercase\expandafter{\romannumeral3}}+{\rm \uppercase\expandafter{\romannumeral4}}\leq (2m+2)^2e^{-s(\sqrt{d}-1)^2}.
\end{equation}
Therefore, \eqref{equ 4.6} follows from \eqref{equ 4.7} and the proof is complete.
\qed

At last, we prove Theorem \ref{theorem 2.1 CLT}.

\proof[Proof of Theorem \ref{theorem 2.1 CLT}]

By Lemmas \ref{lemma 4.1}, \ref{lemma 4.2 martingale part convergence} and \eqref{equ 3.5 extended martingale decomposition formula}, any finite dimensional distribution of
\[
\left\{\frac{1}{\sqrt{N}}\xi_{tN}^F:~0\leq t\leq T\right\}
\]
converges weakly to the corresponding finite dimensional distribution of $\{\sigma_FB_t\}_{0\leq t\leq T}$ as $N\rightarrow+\infty$. Hence, to complete the proof, we only need to show that
\[
\left\{\frac{1}{\sqrt{N}}\xi_{tN}^F:~0\leq t\leq T\right\}_{N\geq 1}
\]
are tight under the uniform topology of $C[0, T]$. A sufficient condition of the above tightness is that there exists $C_1<+\infty$ independent of $s, t, N$ such that
\[
\mathbb{E}_{\nu_p}\left(\left(\frac{1}{\sqrt{N}}\xi_{tN}^F-\frac{1}{\sqrt{N}}\xi_{sN}^F\right)^4\right)\leq C_1(t-s)^2
\]
for any $N\geq 1$ and $0\leq s, t$. According to the invariance of $\nu_p$, to complete the proof, we only need to show that there exists $C_2<+\infty$ independent of $t, N$ such that
\begin{equation}\label{equ 4.8}
\frac{1}{N^2}\int_0^{tN}dt_1\int_0^{t_1}dt_2\int_0^{t_2}dt_3\int_0^{t_3}\mathbb{E}_{\nu_p}\left(\prod_{i=1}^4F(\eta_{t_i})\right)dt_4
\leq C_2t^2
\end{equation}
for any $t\geq 0$ and $N\geq 1$. By \eqref{equ 3.1 graphical representation},
\[
\prod_{i=1}^4F(\eta_{t_i})=\prod_{i=1}^4H(\eta_0, X_{t_i}^{t_i, x_1}, \ldots, X_{t_i}^{t_i, x_m}).
\]
According to the definition of $\nu_p$ and Assumption \eqref{equ basic mean zero assumption}, for
\[
y_1^1,\ldots, y_m^1, y_1^2, \ldots, y_m^2, y_1^3, \ldots, y_m^3, y_1^4, \ldots, y_m^4\in \mathcal{T}_d,
\]
\[
\mathbb{E}_{\nu_P}\left(\prod_{i=1}^4H(\eta_0, y_1^i, \ldots, y_m^i)\right)\neq 0
\]
only if there exist $1\leq i<j\leq 4$ such that
\[
\left\{y_1^i,\ldots, y_m^i\right\}\bigcap \left\{y_1^j, \ldots, y_m^j\right\}\neq \emptyset
\]
and
\[
\left\{y_1^{\hat{i}},\ldots, y_m^{\hat{i}}\right\}\bigcap \left\{y_1^{\hat{j}}, \ldots, y_m^{\hat{j}}\right\}\neq \emptyset,
\]
where $\{\hat{i}, \hat{j}\}=\{1, 2, 3 ,4\}\setminus \{i, j\}$. Hence,
\[
\left|\mathbb{E}_{\nu_p}\left(\prod_{i=1}^4F(\eta_{t_i})\right)\right|\leq K_H^4\left({\rm \uppercase\expandafter{\romannumeral5}}
+{\rm \uppercase\expandafter{\romannumeral6}}+{\rm \uppercase\expandafter{\romannumeral7}}\right),
\]
where
\begin{align*}
{\rm \uppercase\expandafter{\romannumeral5}}
&=\mathbb{P}\Big(\left\{X_{t_1}^{t_1, x_1}, \ldots, X_{t_1}^{t_1, x_m}\right\}\bigcap\left\{X_{t_2}^{t_2, x_1}, \ldots, X_{t_2}^{t_2, x_m}\right\}\neq \emptyset,\\
&\text{\quad\quad\quad}\left\{X_{t_3}^{t_3, x_1}, \ldots, X_{t_3}^{t_3, x_m}\right\}\bigcap\left\{X_{t_4}^{t_4, x_1}, \ldots, X_{t_4}^{t_4, x_m}\right\}\neq \emptyset\Big),
\end{align*}
\begin{align*}
{\rm \uppercase\expandafter{\romannumeral6}}
&=\mathbb{P}\Big(\left\{X_{t_1}^{t_1, x_1}, \ldots, X_{t_1}^{t_1, x_m}\right\}\bigcap\left\{X_{t_3}^{t_3, x_1}, \ldots, X_{t_3}^{t_3, x_m}\right\}\neq \emptyset,\\
&\text{\quad\quad\quad}\left\{X_{t_2}^{t_2, x_1}, \ldots, X_{t_2}^{t_2, x_m}\right\}\bigcap\left\{X_{t_4}^{t_4, x_1}, \ldots, X_{t_4}^{t_4, x_m}\right\}\neq \emptyset\Big)
\end{align*}
and
\begin{align*}
{\rm \uppercase\expandafter{\romannumeral7}}
&=\mathbb{P}\Big(\left\{X_{t_1}^{t_1, x_1}, \ldots, X_{t_1}^{t_1, x_m}\right\}\bigcap\left\{X_{t_4}^{t_4, x_1}, \ldots, X_{t_4}^{t_4, x_m}\right\}\neq \emptyset,\\
&\text{\quad\quad\quad}\left\{X_{t_2}^{t_2, x_1}, \ldots, X_{t_2}^{t_2, x_m}\right\}\bigcap\left\{X_{t_3}^{t_3, x_1}, \ldots, X_{t_3}^{t_3, x_m}\right\}\neq \emptyset\Big).
\end{align*}
According to the definition of $X_s^{t, x}$, for $t>r, x, y\in \mathcal{T}_d$,
\[
X_t^{t, x}=X_r^{r, y}
\]
if and only if $X^{t, x}_{t-r}=y$. Therefore, by \eqref{equ random walk heat estimation},
\[
{\rm \uppercase\expandafter{\romannumeral5}}\leq m^4e^{-(t_1-t_2)(\sqrt{d}-1)^2}e^{-(t_3-t_4)(\sqrt{d}-1)^2},
\]
\[
{\rm \uppercase\expandafter{\romannumeral6}}\leq m^4e^{-(t_1-t_3)(\sqrt{d}-1)^2}e^{-(t_3-t_4)(\sqrt{d}-1)^2}
\]
and
\[
{\rm \uppercase\expandafter{\romannumeral7}}\leq m^4e^{-(t_2-t_3)(\sqrt{d}-1)^2}e^{-(t_3-t_4)(\sqrt{d}-1)^2}.
\]
Consequently,
\begin{align*}
&\frac{1}{N^2}\int_0^{tN}dt_1\int_0^{t_1}dt_2\int_0^{t_2}dt_3\int_0^{t_3}\mathbb{E}_{\nu_p}\left(\prod_{i=1}^4F(\eta_{t_i})\right)dt_4\\
&\leq \frac{3K_H^4}{N^2}m^4(tN)^2\left(\int_0^{+\infty}e^{-u(\sqrt{d}-1)^2}du\right)^2
\end{align*}
and hence \eqref{equ 4.8} holds with
\[
C_2=3m^4K_H^4\left(\int_0^{+\infty}e^{-u(\sqrt{d}-1)^2}du\right)^2.
\]
Since \eqref{equ 4.8} holds, the proof is complete.
\qed

\section{Proof of Theorem \ref{theorem 2.2 mdp}}\label{section five}
In this section, we prove Theorem \ref{theorem 2.2 mdp}. By utilizing the decomposition formula \eqref{equ 3.5 extended martingale decomposition formula}, the proof follows the exponential martingale strategy introduced in \cite{Kipnis1989}, the outline of which we give in Subsection \ref{subsection 5.1}. In Subsection \ref{subsection 5.2}, we prove three lemmas which play the key roles in the exponential martingale strategy. The first lemma shows that $\{\frac{1}{a_t}\xi_t^F\}_{t\geq 0}$ are exponentially tight. The other lemmas shows that the errors in two replacements are super-exponentially small. Proofs of all three lemmas utilize a time version block-analysis introduced in \cite{Xue2025}.

\subsection{Outline of the proof} \label{subsection 5.1}
In this subsection, we give the outline of the proof of Theorem \ref{theorem 2.2 mdp}.  For any $t, s>0$ and $c\in \mathbb{R}$, we define
\[
\mathcal{M}_s^{t, c}=\frac{e^{\frac{ca_t}{t}G_{\sqrt{1/t}}^F(\eta_s)}}{e^{\frac{ca_t}{t}G_{\sqrt{1/t}}^F(\eta_0)}}\exp\left\{-\int_0^s
\frac{\mathcal{L}e^{\frac{ca_t}{t}G_{\sqrt{1/t}}^F(\eta_u)}}{e^{\frac{ca_t}{t}G_{\sqrt{1/t}}^F(\eta_u)}}du\right\},
\]
where $G_\lambda^F$ is defined as in Section \ref{section three}. According to the Feynman-Kac formula, $\{\mathcal{M}_s^{t, c}\}_{s\geq 0}$ is a martingale.  According to the definition of $\mathcal{L}$ and the Taylor's expansion formula up to the second order,
\begin{align*}
&\frac{\mathcal{L}e^{\frac{ca_t}{t}G_{\sqrt{1/t}}^F(\eta_u)}}{e^{\frac{ca_t}{t}G_{\sqrt{1/t}}^F(\eta_u)}}\\
&=\frac{1}{2}\sum_{y\in \mathcal{T}_d}\sum_{z\sim y}\left(\frac{ca_t}{t}\left(G_{\sqrt{1/t}}^F(\eta_u^{y, z})-G_{\sqrt{1/t}}^F(\eta_u)\right)+\frac{c^2a_t^2}{2t^2}\left(G_{\sqrt{1/t}}^F(\eta_u^{y, z})-G_{\sqrt{1/t}}^F(\eta_u)\right)^2\right)\\
&\text{\quad\quad\quad}+o(a_t^2/t^2)\\
&=\frac{ca_t}{t}\mathcal{L}G_{\sqrt{1/t}}^F(\eta_u)+\frac{c^2a_t^2}{4t^2}\sum_{y\in \mathcal{T}_d}\sum_{z\sim y}\left(G_{\sqrt{1/t}}^F(\eta_u^{y, z})-G_{\sqrt{1/t}}^F(\eta_u)\right)^2
+o(a_t^2/t^2).
\end{align*}
Then, by \eqref{equ 3.4 generator on G},
\begin{align*}
&\mathcal{M}_t^{t, c}\\
&=\exp\left\{\frac{a_t^2}{t}\left(\frac{c}{a_t}\xi_t^F-\frac{c^2}{2t}\int_0^t\frac{1}{2}\sum_{y\in \mathcal{T}_d}\sum_{z\sim y}\left(G_{\sqrt{1/t}}^F(\eta_u^{y, z})-G_{\sqrt{1/t}}^F(\eta_u)\right)^2du+c\varepsilon_{1, t}+o(1)\right)\right\},
\end{align*}
where
\[
\varepsilon_{1, t}=\frac{1}{a_t}G_{\sqrt{1/t}}^F(\eta_t)-\frac{1}{a_t}G_{\sqrt{1/t}}^F(\eta_0)-\frac{1}{\sqrt{t}}\int_0^t\frac{1}{a_t}G_{\sqrt{1/t}}^F(\eta_u)du.
\]
To utilize the exponential martingale strategy, we need three lemmas, which are proved in Subsection \ref{subsection 5.2}. The first lemma gives the exponential tightness of $\{\frac{1}{a_t}\xi_t^F\}_{t\geq 0}$.
\begin{lemma}\label{lemma 5.1 exponentially tightness}
\[
\limsup_{M\rightarrow+\infty}\limsup_{t\rightarrow+\infty}\frac{t}{a_t^2}\log \mathbb{P}_{\nu_p}
\left(\left|\frac{1}{a_t}\xi_t^F\right|\geq M\right)=-\infty.
\]
\end{lemma}
The second lemma shows that $\varepsilon_{1, t}$ is super-exponentially small as $t\rightarrow+\infty$.
\begin{lemma}\label{lemma 5.2}
For any $\epsilon>0$,
\[
\limsup_{t\rightarrow+\infty}\frac{t}{a_t^2}\log \mathbb{P}_{\nu_p}\left(|\varepsilon_{1, t}|\geq \epsilon\right)=-\infty.
\]
\end{lemma}
The third lemma shows that, if the term $\sum_{y\in \mathcal{T}_d}\sum_{z\sim y}\left(G_{\sqrt{1/t}}^F(\eta_u^{y, z})-G_{\sqrt{1/t}}^F(\eta_u)\right)^2$ in the expression of $\mathcal{M}_t^{t, c}$ is replaced by its mean under $\nu_p$, then the error is super-exponentially small.
\begin{lemma}\label{lemma 5.3 replacement lemma}
Let
\begin{align*}
&\varepsilon_{2, t}=\frac{1}{t}\int_0^t\sum_{y\in \mathcal{T}_d}\sum_{z\sim y}\Bigg(\left(G_{\sqrt{1/t}}^F(\eta_u^{y, z})-G_{\sqrt{1/t}}^F(\eta_u)\right)^2
\\
&\text{\quad\quad\quad\quad\quad\quad}-\mathbb{E}_{\nu_p}\left(\left(G_{\sqrt{1/t}}^F(\eta_0^{y, z})-G_{\sqrt{1/t}}^F(\eta_0)\right)^2\right)\Bigg)du,
\end{align*}
then, for any $\epsilon>0$,
\[
\limsup_{t\rightarrow+\infty}\frac{t}{a_t^2}\log \mathbb{P}_{\nu_p}\left(|\varepsilon_{2, t}|\geq \epsilon\right)=-\infty.
\]
\end{lemma}
We prove Lemmas \ref{lemma 5.1 exponentially tightness}, \ref{lemma 5.2} and \ref{lemma 5.3 replacement lemma} in Subsection \ref{subsection 5.2}. According to an analysis similar with that leading to \eqref{equ 4.3}, we have
\[
\lim_{t\rightarrow+\infty}\mathbb{E}_{\nu_p}\left(\sum_{y\in \mathcal{T}_d}\sum_{z\sim y}\left(G_{\sqrt{1/t}}^F(\eta_0^{y, z})-G_{\sqrt{1/t}}^F(\eta_0)\right)^2\right)=2\sigma_F^2.
\]
Then, by Lemmas \ref{lemma 5.2} and \ref{lemma 5.3 replacement lemma}, we have
\begin{equation}\label{equ 5.1}
\mathcal{M}_s^{t, c}=\exp\left\{\frac{a_t^2}{t}\left(\frac{c}{a_t}\xi_t^F-\frac{c^2}{2}\sigma_F^2+\varepsilon_{3,t}\right)\right\},
\end{equation}
where
\[
\limsup_{t\rightarrow+\infty}\frac{t}{a_t^2}\log \mathbb{P}_{\nu_p}\left(|\varepsilon_{3, t}|\geq \epsilon\right)=-\infty.
\]
With Lemma \ref{lemma 5.1 exponentially tightness} and Equation \eqref{equ 5.1}, according to the fact that
\[
\frac{u^2}{2\sigma_F^2}=\sup_{c\in \mathbb{R}}\left\{cu-\frac{c^2}{2}\sigma_F^2\right\},
\]
the proof of Theorem \ref{theorem 2.2 mdp} follows from a routine procedure, which is introduced in \cite{Kipnis1989} and utilized in many references such as
\cite{Franco2017, Gao2003, Gao2024, Wang2006, XueZhao2021, Xue2025}. So in this paper we omit the details of how to utilize \eqref{equ 5.1} to prove Theorem \ref{theorem 2.2 mdp}.

\subsection{Proofs of Lemmas \ref{lemma 5.1 exponentially tightness}, \ref{lemma 5.2} and \ref{lemma 5.3 replacement lemma}}\label{subsection 5.2}
We first prove Lemma \ref{lemma 5.1 exponentially tightness}. As a preliminary, for any $t\geq 0$, we define $\vartheta_t: \{0, 1\}^{\mathcal{T}_d}\rightarrow \mathbb{R}$ as
\[
\vartheta_t(\eta)=\mathbb{E}_\eta F(\eta_t)
\]
for any $\eta\in \{0, 1\}^{\mathcal{T}_d}$. We need the following lemma to prove Lemma \ref{lemma 5.1 exponentially tightness}.

\begin{lemma}\label{lemma 5.2.1}
Under $\nu_p$,
\[
\lim_{t\rightarrow+\infty}\vartheta_t=0
\]
in $L^2$.
\end{lemma}

\proof
By \eqref{equ 3.2 mean calculation},
\begin{align*}
&\mathbb{E}_{\nu_p}\left(\vartheta_t^2\right)=\sum_{y_1,\ldots, y_m}\sum_{z_1, \ldots, z_m}\Bigg(\mathbb{E}_{\nu_p}\left(H(\eta_0, y_1, \ldots, y_m)H(\eta_0, z_1, \ldots, y_m)\right)\\
&\text{\quad\quad\quad\quad\quad\quad\quad\quad\quad}\times Q_{t, x_1, \ldots, x_m}^m(y_1, \ldots, y_m)Q_{t, x_1, \ldots, x_m}^m(z_1, \ldots, z_m)\Bigg).
\end{align*}
According to Assumption \eqref{equ basic mean zero assumption},
\[
\mathbb{E}_{\nu_p}\left(H(\eta_0, y_1, \ldots, y_m)H(\eta_0, z_1, \ldots, y_m)\right)=0
\]
when $\{y_1, \ldots, y_m\}\bigcap\{z_1, \ldots, z_m\}=\emptyset$. Hence, by \eqref{equ random walk heat estimation},
\begin{align*}
\mathbb{E}_{\nu_p}\left(\vartheta_t^2\right)
&\leq K_H^2\mathbb{P}\left(\left\{\hat{X}_t^{t, x_1},\ldots, \hat{X}_t^{t, x_m}\right\}\bigcap \left\{X_t^{t, x_1},\ldots, X_t^{t, x_m}\right\}\neq \emptyset\right)\\
&\leq K_H^2m^2e^{-t(\sqrt{d}-1)^2},
\end{align*}
where $K_H, X_t^{t, x}, \hat{X}_t^{t, x}$ are defined as in Section \ref{section four}. Consequently, Lemma \ref{lemma 5.2.1} holds.  \qed

For any $f: \{0, 1\}^{\mathcal{T}_d}\rightarrow [0, +\infty)$, we call $f$ a $\nu_p$-density when
\[
\int_{\{0, 1\}^{\mathcal{T}_d}}f(\eta)\nu_p(d\eta)=1.
\]
By Lemma \ref{lemma 5.2.1}, we have the following lemma.
\begin{lemma}\label{lemma 5.2.2}
For any $\nu_p$-density $f$,
\[
\lim_{t\rightarrow+\infty}\int_{\{0, 1\}^{\mathcal{T}_d}}f(\eta)\vartheta_t(\eta)\nu_p(d\eta)=0.
\]
\end{lemma}

\proof

For any $c>0$, according to the fact that $|\vartheta_t|\leq K_H$, we have
\[
\left|\int_{f\geq c}f(\eta)\vartheta_t(\eta)\nu_p(d\eta)\right|\leq
K_H\int_{f\geq c}f(\eta)\nu_p(d\eta).
\]
Hence, by Lemma \ref{lemma 5.2.1}, Cauchy-Schwartz inequality and the fact that
\[
\left|\int_{f\leq c}\vartheta_t(\eta)f(\eta)\nu_p(d\eta)\right|\leq c\mathbb{E}_{\nu_p}|\vartheta_t|,
\]
we have
\[
\limsup_{t\rightarrow+\infty}\left|\int_{\{0, 1\}^{\mathcal{T}_d}}f(\eta)\vartheta_t(\eta)\nu_p(d\eta)\right|
\leq K_H\int_{f\geq c}f(\eta)\nu_p(d\eta)
\]
for any $c>0$. Since $c$ is arbitrary  and $f$ is a $\nu_p$-density, let $c\rightarrow+\infty$ in the above inequality and then the proof is complete.  \qed

Now we prove Lemma \ref{lemma 5.1 exponentially tightness}.

\proof[Proof of Lemma \ref{lemma 5.1 exponentially tightness}]
We only show that
\[
\limsup_{M\rightarrow+\infty}\limsup_{t\rightarrow+\infty}\frac{t}{a_t^2}\log \mathbb{P}_{\nu_p}\left(\frac{1}{a_t}\xi_t^F\geq M\right)=-\infty,
\]
since
\[
\limsup_{M\rightarrow+\infty}\limsup_{t\rightarrow+\infty}\frac{t}{a_t^2}\log \mathbb{P}_{\nu_p}\left(\frac{1}{a_t}\xi_t^F\leq -M\right)=-\infty
\]
follows from a similar analysis. By Markov inequality,
\begin{equation}\label{equ 5.1 two}
\mathbb{P}_{\nu_p}\left(\frac{1}{a_t}\xi_t^F\geq M\right)\leq e^{-\frac{a_t^2}{t}M}\mathbb{E}_{\nu_p}\exp\left\{\frac{a_t}{t}
\int_0^tF(\eta_s)ds\right\}.
\end{equation}
Since $\nu_p$ is a reversible distribution of $\{\eta_t\}_{t\geq 0}$, by Lemma 7.2 in Appendix 1 of \cite{kipnis+landim99}, we have
\begin{align*}
&\mathbb{E}_{\nu_p}\exp\left\{\frac{a_t}{t}
\int_0^tF(\eta_s)ds\right\}\\
&\leq \exp\left\{t\left(\sup_{f\text{~is a~}\nu_p\text{-density}}\left\{\frac{a_t}{t}\int_{\{0, 1\}^{\mathcal{T}_d}}
f(\eta)F(\eta)\nu_p(d\eta)-\mathcal{D}(\sqrt{f})\right\}\right)\right\},
\end{align*}
where
\[
\mathcal{D}(\sqrt{f})=\frac{1}{4}\sum_{y\in \mathcal{T}_d}\sum_{z\sim y}\int_{\{0, 1\}^{\mathcal{T}_d}}\left(\sqrt{f(\eta^{y, z})}-\sqrt{f(\eta)}\right)^2\nu_p(d\eta),
\]
i.e., the Dirichlet form on $\sqrt{f}$. Then, by \eqref{equ 5.1 two}, to complete the proof we only need to show that
\begin{equation}\label{equ 5.2}
\limsup_{t\rightarrow+\infty}\sup_{f\text{~is a~}\nu_p\text{-density}}\left\{\frac{t}{a_t}\int_{\{0, 1\}^{\mathcal{T}_d}}
f(\eta)F(\eta)\nu_p(d\eta)-\frac{t^2}{a_t^2}\mathcal{D}(\sqrt{f})\right\}<+\infty.
\end{equation}
Now we check \eqref{equ 5.2}. For any $M>0$, by Chapman-Kolmogorov equation,
\begin{align*}
F(\eta)&=\vartheta_0(\eta)=-\int_0^M\frac{d}{dt}\vartheta_t(\eta)dt+\vartheta_M(\eta)\\
&=-\int_0^M\mathcal{L}\vartheta_t(\eta)dt+\vartheta_M(\eta)
\end{align*}
and hence
\begin{align}\label{equ 5.3}
&\int_{\{0, 1\}^{\mathcal{T}_d}}
f(\eta)F(\eta)\nu_p(d\eta)\\
&=-\int_0^M\left(\int_{\{0, 1\}^{\mathcal{T}_d}}
f(\eta)\mathcal{L}\vartheta_t(\eta)\nu_p(d\eta)\right)dt+\int_{\{0, 1\}^{\mathcal{T}_d}}f(\eta)\vartheta_M(\eta)\nu_p(d\eta). \notag
\end{align}
According to the reversibility of $\nu_p$,
\begin{align*}
&\int_{\{0, 1\}^{\mathcal{T}_d}}f(\eta)\mathcal{L}\vartheta_t(\eta)\nu_p(d\eta)\\
&=\frac{1}{2}\int_{\{0, 1\}^{\mathcal{T}_d}}f(\eta)\mathcal{L}\vartheta_t(\eta)+\vartheta_t(\eta)\mathcal{L}f(\eta)
-\mathcal{L}(f\vartheta_t)(\eta)\nu_p(d\eta)\\
&=-\frac{1}{4}\sum_{y\in \mathcal{T}_d}\sum_{z\sim y}\int_{\{0, 1\}^{\mathcal{T}_d}}\left(f(\eta^{y, z})-f(\eta)\right)\left(\vartheta_t(\eta^{y, z})-\vartheta_t(\eta)\right)\nu_p(\eta).
\end{align*}
Therefore,
\begin{align}\label{equ 5.3 two}
&-\int_0^M\left(\int_{\{0, 1\}^{\mathcal{T}_d}}
f(\eta)\mathcal{L}\vartheta_t(\eta)\nu_p(d\eta)\right)dt\\
&=\frac{1}{4}\sum_{y\in \mathcal{T}_d}\sum_{z\sim y}\int_{\{0, 1\}^{\mathcal{T}_d}}\left(f(\eta^{y, z})-f(\eta)\right)\left(A^M(\eta^{y, z})-A^M(\eta)\right)\nu_p(\eta), \notag
\end{align}
where
\[
A^M(\eta)=\int_0^M\vartheta_t(\eta)dt.
\]
According to Cauchy-Schwartz inequality and the fact that
\[
f(\eta^{y, z})-f(\eta)=\left(\sqrt{f(\eta^{y, z})}-\sqrt{f(\eta)}\right)\left(\sqrt{f(\eta^{y, z})}+\sqrt{f(\eta)}\right),
\]
we have
\begin{align}\label{equ 5.4}
&\left|\frac{1}{4}\sum_{y\in \mathcal{T}_d}\sum_{z\sim y}\int_{\{0, 1\}^{\mathcal{T}_d}}\left(f(\eta^{y, z})-f(\eta)\right)\left(A^M(\eta^{y, z})-A^M(\eta)\right)\nu_p(\eta)\right|\\
&\leq \frac{1}{2}\sqrt{\mathcal{D}(\sqrt{f})}\sqrt{\sum_{y\in \mathcal{T}_d}\sum_{z\sim y}\int_{\{0, 1\}^{\mathcal{T}_d}}
\left(\sqrt{f(\eta^{y, z})}+\sqrt{f(\eta)}\right)^2\left(A^M(\eta^{y, z})-A^M(\eta)\right)^2\nu_p(d\eta)}. \notag
\end{align}
By \eqref{equ 3.2 mean calculation},
\begin{align*}
&\sum_{y\in \mathcal{T}_d}\sum_{z\sim y}\left(A^M(\eta^{y, z})-A^M(\eta)\right)^2\\
&=\sum_{y_1,\ldots, y_m}\sum_{z_1, \ldots, z_m}\sum_{y}\sum_{z\sim y}
\int_0^M\int_0^M\Bigg(\Delta_H(\eta, y, z, y_1,\ldots, y_m)\Delta_H(\eta, y, z, z_1, \ldots, z_m)
\\
&\text{\quad\quad}\times Q^m_{t, x_1, \ldots, x_m}(y_1, \ldots, y_m)Q^m_{s, x_1, \ldots, x_m}(z_1, \ldots, z_m)\Bigg)dtds,
\end{align*}
where
\[
\Delta_H(\eta, y, z, y_1,\ldots, y_m)=H(\eta^{y, z}, y_1, \ldots, y_m)-H(\eta, y_1, \ldots, y_m)
\]
for any $y\sim z$ and $y_1, \ldots, y_m$. Note that, when $\{y, z\}\bigcap \{y_1, \ldots, y_m\}=\emptyset$, we have
\[
\Delta_H(\eta, y, z, y_1,\ldots, y_m)=0.
\]
Furthermore, if there exist $y\sim z$ such that
\[
\{y, z\}\bigcap \{y_1, \ldots, y_m\}\neq\emptyset \text{~and~} \{y, z\}\bigcap \{z_1, \ldots, z_m\}\neq \emptyset,
\]
then there exist $1\leq i, j\leq m$ such that $D(y_i, z_j)\leq 1$. Therefore,
\begin{align*}
&\sum_{y_1,\ldots, y_m}\sum_{z_1, \ldots, z_m}\sum_{y}\sum_{z\sim y}
\Bigg(\Delta_H(\eta, y, z, y_1,\ldots, y_m)\Delta_H(\eta, y, z, z_1, \ldots, z_m)
\\
&\text{\quad\quad}\times Q^m_{t, x_1, \ldots, x_m}(y_1, \ldots, y_m)Q^m_{s, x_1, \ldots, x_m}(z_1, \ldots, z_m)\Bigg)\\
&\leq 4K_H^2 \mathbb{P}\left(D\left(X_t^{u, x_i}, \hat{X}_s^{u, x_j}\right)\leq 1\text{~for some~}1\leq i,j\leq m\right),
\end{align*}
where $u\geq t, s$. As we have shown in the proof of Lemma \ref{lemma 4.2 martingale part convergence},
\[
\mathbb{P}\left(D\left(X_t^{u, x_i}, \hat{X}_s^{u, x_j}\right)\leq 1\text{~for some~}1\leq i,j\leq m\right)
\leq m^2(d+2)e^{-\max\{s, t\}(\sqrt{d}-1)^2}.
\]
Therefore, for any $M>0$ and $\eta\in \{0, 1\}^{\mathcal{T}_d}$,
\[
\sum_{y\in \mathcal{T}_d}\sum_{z\sim y}\left(A^M(\eta^{y, z})-A^M(\eta)\right)^2\leq \mathcal{J}_1,
\]
where
\[
\mathcal{J}_1=4K_H^2m^2(d+2)\int_0^{+\infty}\int_0^{+\infty}e^{-\max\{s, t\}(\sqrt{d}-1)^2} ds dt<+\infty.
\]
Then, according to the inequality $(a+b)^2\leq 2a^2+2b^2$ and the fact that $\eta\rightarrow \eta^{y, z}$ is a measure-preserving transformation under $\nu_p$, we have
\begin{align*}
\sum_{y\in \mathcal{T}_d}\sum_{z\sim y}\int_{\{0, 1\}^{\mathcal{T}_d}}
\left(\sqrt{f(\eta^{y, z})}+\sqrt{f(\eta)}\right)^2\left(A^M(\eta^{y, z})-A^M(\eta)\right)^2\nu_p(d\eta)
\leq 4\mathcal{J}_1.
\end{align*}
As a result, by \eqref{equ 5.4}, we have
\[
\left|-\int_0^M\left(\int_{\{0, 1\}^{\mathcal{T}_d}}
f(\eta)\mathcal{L}\vartheta_t(\eta)\nu_p(d\eta)\right)dt\right|\leq \sqrt{\mathcal{D}(\sqrt{f})}\sqrt{\mathcal{J}_1}.
\]
Then, by \eqref{equ 5.3} and Lemma \ref{lemma 5.2.2},
\begin{equation}\label{equ 5.5 two}
\int_{\{0, 1\}^{\mathcal{T}_d}}
f(\eta)F(\eta)\nu_p(d\eta)\leq \sqrt{\mathcal{D}(\sqrt{f})}\sqrt{\mathcal{J}_1}
\end{equation}
for any $\nu_p$-density $f$. Consequently,
\begin{align*}
&\sup_{f\text{~is a~}\nu_p\text{-density}}\left\{\frac{t}{a_t}\int_{\{0, 1\}^{\mathcal{T}_d}}
f(\eta)F(\eta)\nu_p(d\eta)-\frac{t^2}{a_t^2}\mathcal{D}(\sqrt{f})\right\}\\
&\leq \frac{t}{a_t} \sqrt{\mathcal{D}(\sqrt{f})}\sqrt{\mathcal{J}_1}-\frac{t^2}{a_t^2}\mathcal{D}(\sqrt{f})\leq \sup_{u\in \mathbb{R}}\left\{\sqrt{\mathcal{J}_1}u-u^2\right\}=\frac{\mathcal{J}_1}{4}
\end{align*}
and hence \eqref{equ 5.2} holds.
\qed

To prove Lemma \ref{lemma 5.2}, we need the following lemma, which is an analogue of Lemma \ref{lemma 5.2.1}.

\begin{lemma}\label{lemma 5.2.3}
For $u, t>0$, let $\Xi_u^t: \{0, 1\}^{\mathcal{T}_d}\rightarrow \mathbb{R}$ such that
\[
\Xi_u^t(\eta)=\mathbb{E}_{\eta}G_{\sqrt{1/t}}^F(\eta_u)
\]
for any $\eta\in \{0, 1\}^{\mathbb{T}_d}$, then
\[
\lim_{u\rightarrow+\infty}\Xi_u^t=0
\]
in $L^2$ under $\nu_p$.

\end{lemma}

\proof
By \eqref{equ 3.2 mean calculation},
\begin{align*}
\Xi_u^t(\eta)&=\sum_{y_1, \ldots, y_m}\mathbb{E}_\eta H(\eta_u, y_1, \ldots, y_m)\beta_{\sqrt{1/t}, x_1, \ldots, x_m}(y_1, \ldots, y_m)\\
&=\sum_{y_1, \ldots, y_m}\sum_{z_1, \ldots, z_m}H(\eta, z_1, \ldots, z_m)Q^m_{u, y_1, \ldots, y_m}(z_1, \ldots, z_m)\beta_{\sqrt{1/t}, x_1, \ldots, x_m}(y_1, \ldots, y_m).
\end{align*}
Hence,
\begin{align*}
\mathbb{E}_{\nu_p}\left(\left(\Xi_u^t\right)^2\right)
&=\sum_{y_1, \ldots, y_m}\sum_{z_1, \ldots, z_m}
\sum_{w_1, \ldots, w_m}\sum_{v_1, \ldots, v_m}\Bigg(\mathbb{E}_{\nu_p}\big(H(\eta_0, z_1, \ldots, z_m)H(\eta_0, v_1, \ldots, v_m)\big)\\
&\times Q^m_{u, y_1, \ldots, y_m}(z_1, \ldots, z_m)Q^m_{u, w_1, \ldots, w_m}(v_1, \ldots, v_m)\\
&\times\beta_{\sqrt{1/t}, x_1, \ldots, x_m}(y_1, \ldots, y_m)
\beta_{\sqrt{1/t}, x_1, \ldots, x_m}(w_1, \ldots, w_m)\Bigg).
\end{align*}
Note that, by Assumption \eqref{equ basic mean zero assumption},
\[
\mathbb{E}_{\nu_p}\big(H(\eta_0, z_1, \ldots, z_m)H(\eta_0, v_1, \ldots, v_m)\big)=0
\]
when $\{z_1, \ldots, z_m\}\bigcap \{v_1, \ldots, v_m\}=\emptyset$. Therefore,
\begin{align*}
\mathbb{E}_{\nu_p}\left(\left(\Xi_u^t\right)^2\right)
&\leq\sum_{y_1, \ldots, y_m}\sum_{z_1, \ldots, z_m}
\sum_{w_1, \ldots, w_m}\sum_{(v_1, \ldots, v_m)\in \Lambda_{z_1, \ldots, z_m}}\Bigg(K_H^2\\
&\times Q^m_{u, y_1, \ldots, y_m}(z_1, \ldots, z_m)Q^m_{u, w_1, \ldots, w_m}(v_1, \ldots, v_m)\\
&\times\beta_{\sqrt{1/t}, x_1, \ldots, x_m}(y_1, \ldots, y_m)
\beta_{\sqrt{1/t}, x_1, \ldots, x_m}(w_1, \ldots, w_m)\Bigg).
\end{align*}
Note that
\begin{align*}
&\sum_{y_1, \ldots, y_m}\sum_{z_1, \ldots, z_m}
\sum_{w_1, \ldots, w_m}\sum_{(v_1, \ldots, v_m)\in \Lambda_{z_1, \ldots, z_m}}\Bigg(K_H^2\\
&\times Q^m_{u, y_1, \ldots, y_m}(z_1, \ldots, z_m)Q^m_{u, w_1, \ldots, w_m}(v_1, \ldots, v_m)\\
&\times\beta_{\sqrt{1/t}, x_1, \ldots, x_m}(y_1, \ldots, y_m)
\beta_{\sqrt{1/t}, x_1, \ldots, x_m}(w_1, \ldots, w_m)\Bigg)\\
&\leq \sum_{y_1, \ldots, y_m}\sum_{z_1, \ldots, z_m}
\sum_{w_1, \ldots, w_m}\sum_{(v_1, \ldots, v_m)\in \Lambda_{z_1, \ldots, z_m}}\Bigg(K_H^2\\
&\times Q^m_{u, y_1, \ldots, y_m}(z_1, \ldots, z_m)Q^m_{u, w_1, \ldots, w_m}(v_1, \ldots, v_m)\\
&\times\beta_{0, x_1, \ldots, x_m}(y_1, \ldots, y_m)
\beta_{0, x_1, \ldots, x_m}(w_1, \ldots, w_m)\Bigg)\\
&=K_H^2
\int_0^{+\infty}\int_0^{+\infty}\mathbb{P}\left(
\left\{X_{s+u}^{s+u, x_1}, \ldots, X_{s+u}^{s+u, x_m}\right\}\bigcap\left\{
\hat{X}_{t+u}^{t+u, x_1}, \ldots, \hat{X}_{t+u}^{t+u, x_m}\right\}\neq \emptyset\right)dsdt\\
&\leq K_H^2m^2\int_0^{+\infty}\int_0^{+\infty}\exp\left\{
-\max\{s+u, t+u\}(\sqrt{d}-1)^2\right\}dsdt
\end{align*}
and hence
\[
\lim_{u\rightarrow+\infty}\mathbb{E}_{\nu_p}\left(\left(\Xi_u^t\right)^2\right)=0.
\]
\qed

The following lemma is an analogue of Lemma \ref{lemma 5.2.2}.

\begin{lemma}\label{lemma 5.2.4}
For any $\nu_p$-density $f$,
\[
\lim_{u\rightarrow+\infty}\int_{\{0, 1\}^{\mathcal{T}_d}}f(\eta)\Xi_u^t(\eta)\nu_p(d\eta)=0.
\]
\end{lemma}

By Lemma \ref{lemma 5.2.3}, the proof of Lemma \ref{lemma 5.2.4} follows from the same analysis as that in the proof of Lemma \ref{lemma 5.2.2} and hence we omit the details. Now we prove Lemma \ref{lemma 5.2}.

\proof[Proof of Lemma \ref{lemma 5.2}]
Since
\[
\left|\frac{1}{a_t}G_{\sqrt{1/t}}^F(\eta)\right|\leq \frac{K_H}{a_t}\int_0^{+\infty}e^{-\frac{1}{\sqrt{t}}s}ds=\frac{K_H\sqrt{t}}{a_t}\rightarrow 0
\]
as $t\rightarrow+\infty$, we only need to show that
\begin{equation}\label{equ 5.7 zero}
\limsup_{t\rightarrow+\infty}\frac{t}{a_t^2}\log \mathbb{P}_{\nu_p}\left(\left|\frac{1}{\sqrt{t}}\int_0^t\frac{1}{a_t}G_{\sqrt{1/t}}^F(\eta_u)du\right|\geq \epsilon\right)=-\infty.
\end{equation}
Then, according to an analysis similar with that given in the proof of Lemma \ref{lemma 5.1 exponentially tightness}, to prove \eqref{equ 5.7 zero} we only need to show that, for any $\theta>0$,
\begin{equation}\label{equ 5.5}
\limsup_{t\rightarrow+\infty}\sup_{f\text{~is a~}\nu_p\text{-density}}\left\{\frac{\sqrt{t}\theta}{a_t}\int_{\{0, 1\}^{\mathcal{T}_d}}
f(\eta)G_{\sqrt{1/t}}^F(\eta)\nu_p(d\eta)-\frac{t^2}{a_t^2}\mathcal{D}(\sqrt{f})\right\}\leq 0.
\end{equation}
Now we check \eqref{equ 5.5}. For any $M\geq 0$, since
\[
G_{\sqrt{1/t}}^F(\eta)=\Xi_0^t(\eta)=-\int_0^M \frac{d}{du}\Xi_u^t(\eta)du+\Xi_M^t(\eta),
\]
according to an analysis similar with that leading to \eqref{equ 5.3}, we have
\begin{align}\label{equ 5.6}
&\int_{\{0, 1\}^{\mathcal{T}_d}}
f(\eta)G_{\sqrt{1/t}}^F(\eta)\nu_p(d\eta)\\
&=-\int_0^M\left(\int_{\{0, 1\}^{\mathcal{T}_d}}
f(\eta)\mathcal{L}\Xi_u^t(\eta)\nu_p(d\eta)\right)du+\int_{\{0, 1\}^{\mathcal{T}_d}}f(\eta)\Xi_M^t(\eta)\nu_p(d\eta). \notag
\end{align}
Then, according to an analysis similar with that leading to \eqref{equ 5.3 two}, we have
\begin{align*}
&-\int_0^M\left(\int_{\{0, 1\}^{\mathcal{T}_d}}
f(\eta)\mathcal{L}\Xi_u^t(\eta)\nu_p(d\eta)\right)du\\
&=\frac{1}{4}\sum_{y\in \mathcal{T}_d}\sum_{z\sim y}\int_{\{0, 1\}^{\mathcal{T}_d}}\left(f(\eta^{y, z})-f(\eta)\right)\left(\Phi_t^M(\eta^{y, z})-\Phi_t^M(\eta)\right)\nu_p(\eta), \notag
\end{align*}
where
\[
\Phi_t^M(\eta)=\int_0^M\Xi_u^t(\eta)du.
\]
Similar with the argument given in the proof of Lemma \ref{lemma 5.1 exponentially tightness}, we need to calculate
\[
\sum_{y\in \mathcal{T}_d}\sum_{z\sim y}\left(\Phi^M_t(\eta^{y, z})-\Phi^M_t(\eta)\right)^2.
\]
By \eqref{equ 3.2 mean calculation},
\begin{align*}
\Phi^M_t(\eta^{y, z})-\Phi^M_t(\eta)
&=\int_0^M\sum_{y_1, \ldots, y_m}\sum_{z_1, \ldots, z_m}\Big(\Delta_H(\eta, y, z, z_1, \ldots, z_m)\\
&\text{\quad\quad}\times Q^m_{u, y_1, \ldots, y_m}(z_1, \ldots, z_m)\beta_{\sqrt{1/t}, x_1, \ldots, x_m}(y_1, \ldots, y_m)\Big)du
\end{align*}
and hence
\begin{align*}
&\sum_{y\in \mathcal{T}_d}\sum_{z\sim y}\left(\Phi^M_t(\eta^{y, z})-\Phi^M_t(\eta)\right)^2\\
&=\sum_{y_1, \ldots, y_m}\sum_{z_1, \ldots, z_m}\sum_{w_1, \ldots, w_m}\sum_{v_1, \ldots, v_m}\sum_{y}\sum_{z\sim y}\\
&\text{\quad}\Bigg(\int_0^M\int_0^M\Big(\Delta_H(\eta, y, z, z_1, \ldots, z_m)\Delta_H(\eta, y, z, v_1, \ldots, v_m)\\
&\text{\quad\quad\quad}\times Q_{u_1, y_1, \ldots, y_m}^m(z_1, \ldots, z_m)Q_{u_2, w_1, \ldots, w_m}^m(v_1, \ldots, y_m)\\
&\text{\quad\quad\quad}\times \beta_{\sqrt{1/t}, x_1, \ldots, x_m}(y_1, \ldots, y_m)\beta_{\sqrt{1/t}, x_1, \ldots, x_m}(w_1, \ldots, w_m)\Big)du_1 du_2\Bigg).
\end{align*}
As we have pointed out in the proof of Lemma \ref{lemma 5.1 exponentially tightness}, there exist $y\sim z$ such that
\[
\Delta_H(\eta, y, z, z_1, \ldots, z_m)\Delta_H(\eta, y, z, v_1, \ldots, v_m)\neq 0
\]
only if $D(z_i, v_j)\leq 1$ for some $1\leq i,j\leq m$. Therefore,
\begin{align*}
&\sum_{y\in \mathcal{T}_d}\sum_{z\sim y}\left(\Phi^M_t(\eta^{y, z})-\Phi^M_t(\eta)\right)^2\\
&\leq \sum_{y_1, \ldots, y_m}\sum_{z_1, \ldots, z_m}\sum_{w_1, \ldots, w_m}\sum_{(v_1, \ldots, v_m):\atop
D(v_i, z_j)\leq 1\text{~for some~}i, j}\sum_{y}\sum_{z\sim y}\\
&\text{\quad}\Bigg(\int_0^{+\infty}\int_0^{+\infty}\Big(4K_H^2Q_{u_1, y_1, \ldots, y_m}^m(z_1, \ldots, z_m)Q_{u_2, w_1, \ldots, w_m}^m(v_1, \ldots, v_m)\\
&\text{\quad\quad\quad}\times \beta_{0, x_1, \ldots, x_m}(y_1, \ldots, y_m)\beta_{0, x_1, \ldots, x_m}(w_1, \ldots, w_m)\Big)du_1 du_2\Bigg)
\end{align*}
and hence
\begin{align}\label{equ 5.8 two}
\sum_{y\in \mathcal{T}_d}\sum_{z\sim y}\left(\Phi^M_t(\eta^{y, z})-\Phi^M_t(\eta)\right)^2\leq \mathcal{J}_2,
\end{align}
where
\begin{align*}
\mathcal{J}_2&=4K_H^2\int_0^{+\infty}\int_0^{+\infty}\int_0^{+\infty}\int_0^{+\infty}\\
&\text{\quad\quad}\mathbb{P}\left(D\left(X_{u_1+u_3}^{u_1+u_3, x_i}, \hat{X}_{u_2+u_4}^{u_2+u_4, x_j}\right)\leq 1 \text{~for some~}i, j\right)du_1du_2du_3du_4\\
&\leq 4K_H^2m^2(d+2)\int_0^{+\infty}\int_0^{+\infty}\int_0^{+\infty}\int_0^{+\infty}\\
&\text{\quad\quad}\exp\left\{-\max\{u_1+u_3, u_2+u_4\}(\sqrt{d}-1)^2\right\}du_1du_2du_3du_4<+\infty.
\end{align*}
Then, according to an analysis similar with that leading to \eqref{equ 5.5 two}, which utilizes Cauchy-Schwartz inequality, \eqref{equ 5.6} and Lemma \ref{lemma 5.2.4}, we have
\[
\int_{\{0, 1\}^{\mathcal{T}_d}}
f(\eta)G_{\sqrt{1/t}}^F(\eta)\nu_p(d\eta)\leq \sqrt{\mathcal{D}(\sqrt{f})}\sqrt{\mathcal{J}_2}
\]
for any $\nu_p$-density $f$. Consequently,
\begin{align*}
&\sup_{f\text{~is a~}\nu_p\text{-density}}\left\{\frac{\sqrt{t}\theta}{a_t}\int_{\{0, 1\}^{\mathcal{T}_d}}
f(\eta)G_{\sqrt{1/t}}^F(\eta)\nu_p(d\eta)-\frac{t^2}{a_t^2}\mathcal{D}(\sqrt{f})\right\}\\
&\leq \frac{\sqrt{t}\theta}{a_t}\sqrt{\mathcal{D}(\sqrt{f})}\sqrt{\mathcal{J}_2}-\frac{t^2}{a_t^2}\mathcal{D}(\sqrt{f})\\
&\leq \sup_{u\in \mathbb{R}}\left\{\frac{\theta}{\sqrt{t}}u\sqrt{\mathcal{J}_2}-u^2\right\}=\frac{\theta^2\mathcal{J}_2}{4t}
\end{align*}
and hence \eqref{equ 5.5} holds.  \qed

For any $t,u>0$ and $\eta\in \{0, 1\}^{\mathcal{T}_d}$, we denote by $W_t(\eta)$ the term
\[
\sum_{y\in \mathcal{T}_d}\sum_{z\sim y}\left(\left(G_{\sqrt{1/t}}^F(\eta^{y, z})-G_{\sqrt{1/t}}^F(\eta)\right)^2
-\mathbb{E}_{\nu_p}\left(\left(G_{\sqrt{1/t}}^F(\eta_0^{y, z})-G_{\sqrt{1/t}}^F(\eta_0)\right)^2\right)\right)
\]
and by $\Omega_u^t(\eta)$ the expectation $\mathbb{E}_\eta W_t(\eta_u)$. To prove Lemma \ref{lemma 5.3 replacement lemma}, we need the following lemma, which is an analogue of Lemmas \ref{lemma 5.2.1} and \ref{lemma 5.2.3}.
\begin{lemma}\label{lemma 5.2.5}
Under $\nu_p$, for any $t\geq 0$,
\[
\lim_{u\rightarrow+\infty}\Omega_u^t=0
\]
in $L^2$.
\end{lemma}
\proof
According to the definition of $G_{\sqrt{1/t}}^F(\eta)$, we have
\begin{align*}
&\left(G_{\sqrt{1/t}}^F(\eta^{y, z})-G_{\sqrt{1/t}}^F(\eta)\right)^2
-\mathbb{E}_{\nu_p}\left(G_{\sqrt{1/t}}^F(\eta_0^{y, z})-G_{\sqrt{1/t}}^F(\eta_0)\right)^2\\
&=\sum_{y_1, \ldots, y_m}\sum_{z_1, \ldots, z_m}\Big(\Psi(\eta, y, z, y_1,\ldots, y_m, z_1, \ldots, z_m)\beta_{\sqrt{1/t}, x_1, \ldots, x_m}(y_1, \ldots, y_m)\\
&\text{\quad\quad}\times\beta_{\sqrt{1/t}, x_1, \ldots, x_m}(z_1, \ldots, z_m)\Big),
\end{align*}
where
\begin{align*}
&\Psi(\eta, y, z, y_1,\ldots, y_m, z_1, \ldots, z_m)\\
&=\Delta_H(\eta, y, z, y_1, \ldots, y_m)\Delta_H(\eta, y, z, z_1, \ldots, z_m)\\
&\text{\quad\quad\quad}-\mathbb{E}_{\nu_p}\left(\Delta_H(\eta_0, y, z, y_1, \ldots, y_m)\Delta_H(\eta_0, y, z, z_1, \ldots, z_m)\right).
\end{align*}
Hence,
\begin{align}\label{equ 5.9}
\Omega_u^t(\eta)&=\sum_{y_1, \ldots, y_m}\sum_{z_1, \ldots, z_m}\sum_y\sum_{z\sim y}\Big(\mathbb{E}_\eta\Psi(\eta_u, y, z, y_1,\ldots, y_m, z_1, \ldots, z_m)\\
&\text{\quad\quad}\times\beta_{\sqrt{1/t}, x_1, \ldots, x_m}(y_1, \ldots, y_m)\beta_{\sqrt{1/t}, x_1, \ldots, x_m}(z_1, \ldots, z_m)\Big). \notag
\end{align}
By \eqref{equ 3.2 mean calculation},
\begin{align*}
&\mathbb{E}_\eta\Psi(\eta_u, y, z, y_1,\ldots, y_m, z_1, \ldots, z_m)\\
&=\sum_{w}\sum_{v}\sum_{w_1, \ldots, w_m}\sum_{v_1, \ldots, v_m}\Big(
\Psi(\eta, w, v, w_1, \ldots, w_m, v_1, \ldots, v_m)\\
&\text{\quad\quad}\times Q^{2m+2}_{u, y, z, y_1, \ldots, y_m, z_1, \ldots, z_m}(w, v, w_1, \ldots, w_m, v_1, \ldots, v_m)\Big).
\end{align*}
Note that
\[
\mathbb{E}_{\nu_p}\Psi(\eta_0, w, v, w_1, \ldots, w_m, v_1, \ldots, v_m)=0
\]
for any $w, v, w_1, \ldots, w_m, v_1, \ldots, v_m$. Hence, according to the analysis we have utilized several times,
\begin{align*}
&\mathbb{E}_{\nu_p}\left(\mathbb{E}_{\eta_0}\Psi(\eta_u, y, z, y_1,\ldots, y_m, z_1, \ldots, z_m)\mathbb{E}_{\eta_0}\Psi(\eta_u, \tilde{y}, \tilde{z}, \tilde{y}_1,\ldots, \tilde{y}_m, \tilde{z}_1, \ldots, \tilde{z}_m)\right)\\
&\leq 4K_H^4 \mathbb{P}\left(\left\{X_u^{u, y}, X_u^{u, z}, X_u^{u, y_1}, \ldots, X_u^{u, z_m}\right\}\bigcap \left\{\hat{X}_u^{u, \tilde{y}}, \hat{X}_u^{u, \tilde{z}}, \hat{X}_u^{u, \tilde{y}_1}, \ldots, \hat{X}_u^{u, \tilde{z}_m}\right\}\neq \emptyset\right)\\
&\leq  4K_H^4(2m+2)^2e^{-u(\sqrt{d}-1)^2}
\end{align*}
for any $y\sim z, \tilde{y}\sim \tilde{z}$, $y_1, \ldots, y_m, \ldots, \tilde{z}_1, \ldots, \tilde{z}_m\in \mathcal{T}_d$. According to the definition of $\Delta_H$, as in the proofs of Lemmas \ref{lemma 5.2.1} and \ref{lemma 5.2.3}, in \eqref{equ 5.9} the second sum is still over $(z_1, \ldots, z_m)$ such that
$D(z_i, y_j)\leq 1$ for some $i, j$. As a result,
\begin{align*}
&\mathbb{E}_{\nu_p}\left((\Omega_u^t)^2\right)\\
&\leq \left(\int_0^{+\infty}\int_0^{+\infty}\mathbb{P}\left(
D\left(X_{s+r}^{s+r, x_i}, \hat{X}_{t+r}^{t+r, x_j}\right)\leq 1\text{~for some~}i, j\right)
dsdt\right)^2\\
&\text{\quad\quad}\times  4K_H^4(2m+2)^2e^{-u(\sqrt{d}-1)^2}\\
&\leq 4K_H^4(2m+2)^2\left(m^2(d+2)\int_0^{+\infty}\int_0^{+\infty}\exp\left\{
-\max\{s+r, t+r\}(\sqrt{d}-1)^2\right\}dsdt\right)^2\\
&\text{\quad\quad}\times e^{-u(\sqrt{d}-1)^2}
\end{align*}
and hence Lemma \ref{lemma 5.2.5} holds.

\qed

The following lemma is an analogue of Lemmas \ref{lemma 5.2.2} and \ref{lemma 5.2.4}.
\begin{lemma}\label{lemma 5.2.6}
For any $\nu_p$-density $f$,
\[
\lim_{u\rightarrow+\infty}\int_{\{0, 1\}^{\mathcal{T}_d}}f(\eta)\Omega_u^t(\eta)\nu_p(d\eta)=0.
\]
\end{lemma}
By Lemma \ref{lemma 5.2.5}, the proof of Lemma \ref{lemma 5.2.6} follows from an argument similar with that in the proof of Lemma \ref{lemma 5.2.2}, the detail of which we omit here.

At last, we prove Lemma \ref{lemma 5.3 replacement lemma}.

\proof[Proof of Lemma \ref{lemma 5.3 replacement lemma}]
According to an argument similar with those in proofs of Lemmas \ref{lemma 5.1 exponentially tightness} and \ref{lemma 5.2}, to complete the proof we only need to show that, for any $\theta>0$,
\begin{equation}\label{equ 5.10}
\limsup_{t\rightarrow+\infty}
\sup_{f\text{~is a~}\nu_p\text{-density}}\left\{\theta\int_{\{0, 1\}^{\mathcal{T}_d}}f(\eta)W_t(\eta)\nu_p(d\eta)-\frac{t^2}{a_t^2}\mathcal{D}(\sqrt{f})\right\}\leq 0.
\end{equation}
Now we check \eqref{equ 5.10}. According to the fact that, for any $M>0$,
\[
W_t(\eta)=\Omega_0^t(\eta)=-\int_0^M\frac{d}{du}\Omega_u^t(\eta)du+\Omega_M^t(\eta),
\]
we have
\begin{align}\label{equ 5.11}
&\int_{\{0, 1\}^{\mathcal{T}_d}}
f(\eta)W_t(\eta)\nu_p(d\eta)\\
&=-\int_0^M\left(\int_{\{0, 1\}^{\mathcal{T}_d}}
f(\eta)\mathcal{L}\Omega_u^t(\eta)\nu_p(d\eta)\right)du+\int_{\{0, 1\}^{\mathcal{T}_d}}f(\eta)\Omega_M^t(\eta)\nu_p(d\eta). \notag
\end{align}
According to an analysis similar with that leading to \eqref{equ 5.3 two}, we have
\begin{align*}
&-\int_0^M\left(\int_{\{0, 1\}^{\mathcal{T}_d}}
f(\eta)\mathcal{L}\Omega_u^t(\eta)\nu_p(d\eta)\right)dt\\
&=\frac{1}{4}\sum_{y\in \mathcal{T}_d}\sum_{z\sim y}\int_{\{0, 1\}^{\mathcal{T}_d}}\left(f(\eta^{y, z})-f(\eta)\right)\left(\Upsilon_t^M(\eta^{y, z})-\Upsilon^M_t(\eta)\right)\nu_p(\eta), \notag
\end{align*}
where
\[
\Upsilon_t^M(\eta)=\int_0^M\Omega_u^t(\eta)du.
\]
As in proofs of Lemmas \ref{lemma 5.1 exponentially tightness} and \ref{lemma 5.2}, we need to calculate
\[
\sum_{y\in \mathcal{T}_d}\sum_{z\sim y}\left(\Upsilon_t^M(\eta^{y, z})-\Upsilon^M_t(\eta)\right)^2.
\]
By \eqref{equ 5.9},
\begin{align}\label{equ 5.12}
&\Omega_u^t(\eta^{y, z})-\Omega_u^t(\eta)\notag\\
&=\sum_{y_1, \ldots, y_m}\sum_{z_1, \ldots, z_m}\sum_{\hat{y}}\sum_{\hat{z}\sim \hat{y}}\sum_{w}\sum_{v}\sum_{w_1, \ldots, w_m}\sum_{v_1, \ldots, v_m}\Big(\notag\\
&\text{\quad}\Delta_\Psi(\eta, y, z, w, v, w_1,\ldots, w_m, v_1, \ldots, v_m)Q^{2m+2}_{u, \hat{y}, \hat{z}, y_1, \ldots, y_m, z_1, \ldots, z_m}
(w, v, w_1, \ldots, w_m, v_1, \ldots, v_m)\notag\\
&\times \beta_{\sqrt{1/t}, x_1, \ldots, x_m}(y_1, \ldots, y_m)\beta_{\sqrt{1/t}, x_1, \ldots, x_m}(z_1, \ldots, z_m)\Big),
\end{align}
where
\begin{align*}
&\Delta_\Psi(\eta, y, z, w, v, w_1,\ldots, w_m, v_1, \ldots, v_m)\\
&=\Psi(\eta^{y, z}, w, v, w_1, \ldots, w_m, v_1, \ldots, v_m)-
\Psi(\eta, w, v, w_1, \ldots, w_m, v_1, \ldots, v_m).
\end{align*}
According to the definition of $\Delta_\Psi$, for $y\sim z$,
\[
\Delta_\Psi(\eta, y, z, w, v, w_1,\ldots, w_m, v_1, \ldots, v_m)\neq 0
\]
only if
\[
\{y, z\}\bigcap\left\{w, v, w_1,\ldots, w_m, v_1, \ldots, v_m\right\}\neq \emptyset.
\]
Hence, for $y\sim z$ and $w_0, w_1,\ldots, w_m, w_{m+1}, \ldots, w_{2m+1}, \tilde{w}_0, \tilde{w}_1, \ldots, \tilde{w}_{2m+1}\in \mathcal{T}_d$,
\[
\Delta_\Psi(\eta, y, z, w_0, w_1,\ldots, w_m, w_{m+1}, \ldots, w_{2m+1})
\Delta_\Psi(\eta, y, z, \tilde{w}_0, \tilde{w}_1,\ldots, \tilde{w}_{2m+1})\neq 0
\]
only if $D(w_i, \tilde{w}_j)\leq 1$ for some $0\leq i, j\leq 2m+1$. As in the proof of Lemma \ref{lemma 5.2}, the second sum in \eqref{equ 5.12} is over $(z_1, \ldots, z_m)$ such that $D(z_i, y_j)\leq 1$ for some $1\leq i\leq j\leq m$. Then, according to an analysis similar with that leading to \eqref{equ 5.8 two}, we have
\[
\sum_{y\in \mathcal{T}_d}\sum_{z\sim y}\left(\Upsilon_t^M(\eta^{y, z})-\Upsilon^M_t(\eta)\right)^2\leq \mathcal{J}_3=4K_H^4\mathcal{J}_4\mathcal{J}_5,
\]
where
\begin{align*}
\mathcal{J}_4&=\left(\int_0^{+\infty}\int_0^{+\infty}\mathbb{P}\left(
D\left(X_{s+r}^{s+r, x_i}, \hat{X}_{t+r}^{t+r, x_j}\right)\leq 1\text{~for some~}i, j\right)
dsdt\right)^2\\
&\leq\left(m^2(d+2)\int_0^{+\infty}\int_0^{+\infty}\exp\left\{
-\max\{s+r, t+r\}(\sqrt{d}-1)^2\right\}dsdt\right)^2
\end{align*}
and
\begin{align*}
\mathcal{J}_5=(2m+2)^2(d+2)\int_0^{+\infty}\int_0^{+\infty}e^{-\max\{u_1, u_2\}(\sqrt{d}-1)^2}du_1du_2,
\end{align*}
which is the upper bound of
\begin{align*}
\int_0^{+\infty}\int_0^{+\infty}
\mathbb{P}\left(D\left(X_{u_1}^{u_1, y_i}, \hat{X}_{u_2}^{u_2, \tilde{y}_j}\right)\leq 1\text{~for some~}0\leq i, j\leq 2m+1\right)du_1du_2
\end{align*}
for any $y_0, y_1,\ldots, y_m, y_{m+1}, \ldots, y_{2m+1}, \tilde{y}_0, \tilde{y}_1, \ldots, \tilde{y}_{2m+1}\in \mathcal{T}_d$. Then, according to an analysis similar with that leading to \eqref{equ 5.5 two}, which utilizes Cauchy-Schwartz inequality, \eqref{equ 5.11} and Lemma \ref{lemma 5.2.6}, we have
\[
\int_{\{0, 1\}^{\mathcal{T}_d}}f(\eta)W_t(\eta)\nu_p(d\eta)\leq \sqrt{\mathcal{D}(\sqrt{f})}\sqrt{\mathcal{J}_3}.
\]
Consequently,
\begin{align*}
&\sup_{f\text{~is a~}\nu_p\text{-density}}\left\{\theta\int_{\{0, 1\}^{\mathcal{T}_d}}
f(\eta)W_t(\eta)\nu_p(d\eta)-\frac{t^2}{a_t^2}\mathcal{D}(\sqrt{f})\right\}\\
&\leq \theta\sqrt{\mathcal{D}(\sqrt{f})}\sqrt{\mathcal{J}_3}-\frac{t^2}{a_t^2}\mathcal{D}(\sqrt{f})\\
&\leq \sup_{u\in \mathbb{R}}\left\{\theta u\sqrt{\mathcal{J}_3}-\frac{t^2}{a_t^2}u^2\right\}=\frac{a_t^2\theta^2\mathcal{J}_3}{4t^2}
\end{align*}
and hence \eqref{equ 5.10} holds according to the assumption that $a_t/t\rightarrow 0$.  \qed

\section{Proof of Theorem \ref{theorem 2.3 sample path mdp}}\label{section six}
In this section, we give the outline of the proof of Theorem \ref{theorem 2.3 sample path mdp}. Our proof still follows the exponential martingale strategy. For any $g\in C^1[0, T]$, $N\geq 1$ and $0\leq s\leq T$, we define
\[
\widetilde{M}_s^{N, g}=\frac{e^{\frac{b_N}{N}g_sG_{\sqrt{1/N}}^F(\eta_{sN})}}{e^{\frac{b_N}{N}g_0G_{\sqrt{1/N}}^F(\eta_0)}}\exp\left\{-\int_0^s
\frac{\left(N\mathcal{L}+\partial_u\right)e^{\frac{b_N}{N}g_uG_{\sqrt{1/N}}^F(\eta_{Nu})}}{e^{\frac{b_N}{N}g_uG_{\sqrt{1/N}}^F(\eta_{Nu})}}du\right\},
\]
then $\{\widetilde{M}_s^{N, g}\}_{0\leq s\leq T}$ is a martingale according to the Feynman-Kac formula. As we have shown in Section \ref{section five},
\[
G_{\sqrt{1/N}}^F(\eta)=O(\sqrt{N})
\]
and hence, according to the assumption $\sqrt{N\log N}/b_N\rightarrow 0$,
\begin{align}\label{equ 6.1}
\int_0^T
\frac{\partial_ue^{\frac{b_N}{N}g_uG_{\sqrt{1/N}}^F(\eta_{Nu})}}{e^{\frac{b_N}{N}g_uG_{\sqrt{1/N}}^F(\eta_{Nu})}}du
&=\int_0^T \frac{b_N}{N}g_u^\prime G_{\sqrt{1/N}}^F(\eta_{Nu})du \notag\\
&=\frac{b_N^2}{N}O(\frac{\sqrt{N}}{b_N})=\frac{b_N^2}{N}o(1).
\end{align}
Similarly, we have
\begin{equation}\label{equ 6.2}
\frac{e^{\frac{b_N}{N}g_TG_{\sqrt{1/N}}^F(\eta_{TN})}}{e^{\frac{b_N}{N}g_0G_{\sqrt{1/N}}^F(\eta_0)}}=\exp\left\{\frac{b_N^2}{N}o(1)\right\}.
\end{equation}
According to the definition of $\mathcal{L}$ and the Taylor's expansion formula up to the second order,
\begin{align}\label{equ 6.3}
&\int_0^T
\frac{N\mathcal{L}e^{\frac{b_N}{N}g_uG_{\sqrt{1/N}}^F(\eta_{Nu})}}{e^{\frac{b_N}{N}g_uG_{\sqrt{1/N}}^F(\eta_{Nu})}}du\notag\\
&=\frac{1}{2}\int_0^{T}N\sum_{y\in \mathcal{T}_d}\sum_{z\sim y}\left(e^{\frac{b_N}{N}g_u
\left(G_{\sqrt{1/N}}^F(\eta_{Nu}^{y, z})-G_{\sqrt{1/N}}^F(\eta_{Nu})\right)}-1\right)du\notag\\
&=N\int_0^Tg_u\frac{b_N}{N}\mathcal{L}G_{\sqrt{1/N}}^F(\eta_{Nu})du
+\frac{b_N^2}{4N}\int_0^Tg_u^2\sum_{y\in \mathcal{T}_d}\sum_{z\sim y}\left(G_{\sqrt{1/N}}^F(\eta_{Nu}^{y, z})-G_{\sqrt{1/N}}^F(\eta_{Nu})\right)^2du\notag\\
&\text{\quad}+\frac{b_N^2}{N}o(1).
\end{align}
By \eqref{equ 3.4 generator on G},
\begin{align*}
&N\int_0^Tg_u\frac{b_N}{N}\mathcal{L}G_{\sqrt{1/N}}^F(\eta_{Nu})du\\
&=b_N\int_0^Tg_u\left(-F(\eta_{Nu})+\frac{1}{\sqrt{N}}G_{\sqrt{1/N}}^F(\eta_{Nu})\right)du\\
&=\frac{b_N^2}{N}\left(-\int_0^Tg_u\frac{\partial}{\partial_u}\left(\frac{1}{b_N}\xi_{Nu}^F\right)du+\frac{1}{b_N\sqrt{N}}\int_0^{TN}
g_{u/N}G^F_{\sqrt{1/N}}(\eta_{u})du\right).
\end{align*}
According to an analysis similar with that leading to \eqref{equ 5.7 zero}, we have
\begin{equation}\label{equ 6.4}
\limsup_{N\rightarrow+\infty}\frac{N}{b_N^2}\log \mathbb{P}_{\nu_p}\left(\left|\frac{1}{b_N\sqrt{N}}\int_0^{TN}
g_{u/N}G^F_{\sqrt{1/N}}(\eta_{u})du\right|>\epsilon\right)=-\infty
\end{equation}
for any $\epsilon>0$. Let
\begin{align*}
&\varepsilon_{4,N}=
\int_0^Tg_u^2\Bigg(\sum_{y\in \mathcal{T}_d}\sum_{z\sim y}\Bigg(\left(G_{\sqrt{1/N}}^F(\eta_{Nu}^{y, z})-G_{\sqrt{1/N}}^F(\eta_{Nu})\right)^2
\\
&\text{\quad\quad\quad\quad\quad\quad}-\mathbb{E}_{\nu_p}\left(\left(G_{\sqrt{1/N}}^F(\eta_{0}^{y, z})-G_{\sqrt{1/N}}^F(\eta_{0})\right)^2\right)\Bigg)\Bigg)du\\
&=\frac{1}{N}\int_0^{TN}g_{u/N}^2\Bigg(\sum_{y\in \mathcal{T}_d}\sum_{z\sim y}\Bigg(\left(G_{\sqrt{1/N}}^F(\eta_{u}^{y, z})-G_{\sqrt{1/N}}^F(\eta_{u})\right)^2
\\
&\text{\quad\quad\quad\quad\quad\quad}-\mathbb{E}_{\nu_p}\left(\left(G_{\sqrt{1/N}}^F(\eta_{0}^{y, z})-G_{\sqrt{1/N}}^F(\eta_{0})\right)^2\right)\Bigg)\Bigg)du,
\end{align*}
then according to an analysis similar with that leading to Lemma \ref{lemma 5.3 replacement lemma}, we have
\begin{equation}\label{equ 6.5}
\limsup_{N\rightarrow+\infty}\frac{N}{b_N^2}\log \mathbb{P}_{\nu_p}\left(\left|\varepsilon_{4, N}\right|>\epsilon\right)=-\infty
\end{equation}
for any $\epsilon>0$. According to an analysis similar with that leading to \eqref{equ 4.3}, we have
\[
\lim_{N\rightarrow+\infty}\mathbb{E}_{\nu_p}\left(\sum_{y\in \mathcal{T}_d}\sum_{z\sim y}\left(G_{\sqrt{1/N}}^F(\eta_0^{y, z})-G_{\sqrt{1/N}}^F(\eta_0)\right)^2\right)=2\sigma_F^2.
\]
Consequently, by \eqref{equ 6.1}-\eqref{equ 6.5} and the integration-by-parts formula,
\begin{align}\label{equ 6.6}
\widetilde{M}_T^{N, g}&=\exp\left\{
\frac{b_N^2}{N}\left(\int_0^Tg_u\frac{\partial}{du}\left(\frac{1}{b_N}\xi_{Nu}^F\right)du-\frac{1}{2}\int_0^Tg_u^2\sigma_F^2du+\varepsilon_{5, N}\right)\right\}\notag\\
&=\exp\left\{
\frac{b_N^2}{N}\left(\frac{1}{b_N}\xi_{NT}^Fg_T-\int_0^Tg^\prime_u\left(\frac{1}{b_N}\xi_{Nu}^F\right)du-\frac{1}{2}\int_0^Tg_u^2\sigma_F^2du+\varepsilon_{5, N}\right)\right\},
\end{align}
where
\[
\limsup_{N\rightarrow+\infty}\frac{N}{b_N^2}\log \mathbb{P}_{\nu_p}\left(\left|\varepsilon_{5, N}\right|>\epsilon\right)=-\infty
\]
for any $\epsilon>0$. By \eqref{equ equivalent definition of I} and \eqref{equ 6.6}, the proof of \eqref{equ sample path MDP upper bound} for all compact $\mathcal{C}$ and the proof of \eqref{equ sample path MDP lower bound} for all open $\mathcal{O}$ follow from a routine procedure introduced in \cite{Kipnis1989}, which we omit in this paper. At last, to complete the proof of Theorem \ref{theorem 2.3 sample path mdp}, we only need to show that
\[
\left\{\frac{1}{b_N}\xi_{tN}^F:~0\leq t\leq T\right\}_{N\geq 1}
\]
are exponentially tight. According to the criterion given in \cite{Puhalskii1994}, we need the following lemma.

\begin{lemma}\label{lemma 6.1}
We have
\begin{equation}\label{equ 6.7}
\limsup_{M\rightarrow+\infty}\limsup_{N\rightarrow+\infty}\frac{N}{b_N^2}\log \mathbb{P}_{\nu_p}\left(\sup_{0\leq t\leq T}\left|\frac{1}{b_N}\xi_{tN}^F\right|>M\right)=-\infty
\end{equation}
and
\begin{equation}\label{equ 6.7 two}
\limsup_{\delta\rightarrow 0}\limsup_{N\rightarrow+\infty}\frac{N}{b_N^2}
\log \mathbb{P}_{\nu_p}\left(\sup_{0\leq s\leq t\leq T, \atop
t-s\leq \delta}\left|\frac{1}{b_N}\xi_{tN}^F-\frac{1}{b_N}\xi_{sN}^F\right|\geq \epsilon\right)=-\infty
\end{equation}
for any $\epsilon>0$.
\end{lemma}

\proof
We first check \eqref{equ 6.7}. We only show that
\begin{equation}\label{equ 6.8}
\limsup_{M\rightarrow+\infty}\limsup_{N\rightarrow+\infty}\frac{N}{b_N^2}\log \mathbb{P}_{\nu_p}\left(\sup_{0\leq t\leq T}\frac{1}{b_N}\xi_{tN}^F>M\right)=-\infty
\end{equation}
since
\[
\limsup_{M\rightarrow+\infty}\limsup_{N\rightarrow+\infty}\frac{N}{b_N^2}\log \mathbb{P}_{\nu_p}\left(\inf_{0\leq t\leq T}\frac{1}{b_N}\xi_{tN}^F<-M\right)=-\infty
\]
follows from a similar analysis. For each integer $0\leq i\leq N^3-1$ and $s\in \left[\frac{Ti}{N^2}, \frac{T(i+1)}{N^2}\right)$,
\[
\left|\frac{1}{b_N}\xi_{s}^F-\frac{1}{b_N}\xi_{\frac{iT}{N^2}}^F\right|\leq \frac{1}{b_N}\frac{T}{N^2}K_H\leq 1
\]
for sufficiently large $N$. Hence, to prove \eqref{equ 6.8}, we only need to show that
\begin{equation}\label{equ 6.9}
\limsup_{M\rightarrow+\infty}\limsup_{N\rightarrow+\infty}\frac{N}{b_N^2}\log \left(\sum_{i=0}^{N^3-1}\mathbb{P}_{\nu_p}\left(\frac{1}{b_N}\xi_{\frac{Ti}{N^2}}^F>M\right)\right)=-\infty.
\end{equation}
By Lemma 7.2 in Appendix 1 of \cite{kipnis+landim99} and Markov inequality,
\begin{align*}
&\mathbb{P}_{\nu_p}\left(\frac{1}{b_N}\xi_{\frac{Ti}{N^2}}^F>M\right)\\
&\leq \exp\left\{-\frac{b_N^2M}{N}\right\}
\exp\left\{\frac{Ti}{N^2}\sup_{f\text{~is a~}\nu_p\text{-density}}\left\{\frac{b_N}{N}\int_{\{0, 1\}^{\mathcal{T}_d}}
f(\eta)F(\eta)\nu_p(d\eta)-\mathcal{D}(\sqrt{f})\right\}\right\}.
\end{align*}
Then, by \eqref{equ 5.5 two},
\begin{align*}
\mathbb{P}_{\nu_p}\left(\frac{1}{b_N}\xi_{\frac{Ti}{N^2}}^F>M\right)
&\leq \exp\left\{-\frac{b_N^2M}{N}\right\}\exp\left\{\frac{Ti}{N^2}\sup_{u\in \mathbb{R}}\left\{\frac{b_N}{N}u\sqrt{\mathcal{J}_1}-u^2\right\}\right\}\\
&=\exp\left\{-\frac{b_N^2M}{N}\right\}\exp\left\{\frac{Ti}{N^2}\frac{b_N^2\mathcal{J}_1}{4N^2}\right\}\\
&\leq \exp\left\{-\frac{b_N^2M}{N}\right\}\exp\left\{\frac{b_N^2}{4N}T\mathcal{J}_1\right\}
\end{align*}
for all $0\leq i\leq N^3-1$. As a result,
\begin{equation}\label{equ 6.11 zero}
\frac{N}{b_N^2}\log \left(\sum_{i=0}^{N^3-1}\mathbb{P}_{\nu_p}\left(\frac{1}{b_N}\xi_{\frac{Ti}{N^2}}^F>M\right)\right)
\leq \frac{3N\log N}{b_N^2}-M+T\mathcal{J}_1
\end{equation}
and hence \eqref{equ 6.9} holds according to the assumption $\frac{N\log N}{b_N^2}\rightarrow 0$.

Now we check \eqref{equ 6.7 two}. We only show that
\begin{equation}\label{equ 6.11}
\limsup_{\delta\rightarrow 0}\limsup_{N\rightarrow+\infty}\frac{N}{b_N^2}
\log \mathbb{P}_{\nu_p}\left(\sup_{0\leq s\leq t\leq T, \atop
t-s\leq \delta}\left(\frac{1}{b_N}\xi_{tN}^F-\frac{1}{b_N}\xi_{sN}^F\right)\geq \epsilon\right)=-\infty
\end{equation}
since
\[
\limsup_{\delta\rightarrow 0}\limsup_{N\rightarrow+\infty}\frac{N}{b_N^2}
\log \mathbb{P}_{\nu_p}\left(\inf_{0\leq s\leq t\leq T, \atop
t-s\leq \delta}\left(\frac{1}{b_N}\xi_{tN}^F-\frac{1}{b_N}\xi_{sN}^F\right)\leq -\epsilon\right)=-\infty
\]
follows from a similar analysis. According to the reversibility of $\nu_p$ and a $3$-epsilon argument, to prove \eqref{equ 6.11} we only need to show that
\begin{equation}\label{equ 6.12}
\limsup_{\delta\rightarrow 0}\limsup_{N\rightarrow+\infty}\frac{N}{b_N^2}
\log \mathbb{P}_{\nu_p}\left(\sup_{0\leq t\leq \delta}\frac{1}{b_N}\xi_{tN}^F\geq \epsilon\right)=-\infty.
\end{equation}
For any integer $0\leq i\leq N^3-1$ and $s\in [\frac{\delta i}{N^2}, \frac{\delta(i+1)}{N^2})$,
\[
\left|\frac{1}{b_N}\xi_{s}^F-\frac{1}{b_N}\xi^F_{\frac{\delta i}{N^2}}\right|
\leq \frac{\delta}{N^2}\frac{1}{b_N}K_H\leq \frac{\epsilon}{2}
\]
for $\delta\leq 1$ and sufficiently large $N$. Hence, to prove \eqref{equ 6.12}, we only need to show that
\begin{equation}\label{equ 6.13}
\limsup_{\delta\rightarrow 0}\limsup_{N\rightarrow+\infty}\frac{N}{b_N^2}
\log \left(\sum_{i=0}^{N^3-1}\mathbb{P}_{\nu_p}\left(\frac{1}{b_N}\xi_{\frac{\delta i}{N^2}}^F\geq \epsilon\right)\right)=-\infty.
\end{equation}
According to an analysis similar with that leading to \eqref{equ 6.11 zero}, we have
\begin{align*}
\frac{N}{b_N^2}\log \left(\sum_{i=0}^{N^3-1}\mathbb{P}_{\nu_p}\left(\frac{1}{b_N}\xi_{\frac{\delta i}{N^2}}^F\geq \epsilon\right)\right)
&=\frac{N}{b_N^2}\log \left(\sum_{i=0}^{N^3-1}\mathbb{P}_{\nu_p}\left(\frac{\theta}{b_N}\xi_{\frac{\delta i}{N^2}}^F\geq \theta\epsilon\right)\right)\\
&\leq \frac{3N\log N}{b_N^2}-\theta\epsilon+\delta\theta^2\mathcal{J}_1
\end{align*}
for any $\theta>0$. Hence,
\[
\limsup_{\delta\rightarrow 0}\limsup_{N\rightarrow+\infty}\frac{N}{b_N^2}
\log \left(\sum_{i=0}^{N^3-1}\mathbb{P}_{\nu_p}\left(\frac{1}{b_N}\xi_{\frac{\delta i}{N^2}}^F\geq \epsilon\right)\right)\leq -\theta \epsilon.
\]
Since $\theta$ is arbitrary, let $\theta\rightarrow+\infty$ in the above inequality and \eqref{equ 6.13} holds.
\qed

\quad

\textbf{Acknowledgments.}
The author is grateful to Dr. Linjie Zhao for useful comments.
The author is grateful to financial
supports from the National Natural Science Foundation of China with grant number 12371142 and the Fundamental Research Funds for the Central Universities with grant number 2022JBMC039.

{}
\end{document}